\documentclass[11pt,twoside]{article}
\usepackage{amsfonts}
\begin{document}
\centerline{\bf Mathematical Structures Defined by
Identities}\bigskip \centerline{Constantin M. Petridi,}
\centerline{11 Apollonos Street, 151 24 Maroussi (Athens),
Greece, Fax+30-1-6129250} \centerline{In part - collaboration
with} \centerline{Peter B, Krikelis, Mathematics Dep., Athens
University,} \centerline{e-mail pkrikel@cc.uoa.gr}
\par

\begin{abstract} We propound the thesis that there is a ''limitation'' to
the number of possible structures which are axiomatically endowed
with identities involving operations. In the case of algebras
with a binary operation satisfying a formally irreducible (to be
explained) identity between two $n$-iterates of the operation, it
is established that the frequency of such algebras goes to zero
like $e^{-n/16}$ as $n\to\infty$. This is proved by a suitable
ordering and labeling of the expressions (words) of the
corresponding free algebra and the formation of a series of
tableaux whose entries are the labels. The tableaux reveal
surprising symmetry properties, stated in terms of the Catalan
numbers ${\frac{1}{n+1}}{2n \choose n}$ and their partitions. As
a result of the defining identity and the tableaux all iterates
of order higher than $n$ fall into equivalence classes of
semantically equal ones. Classnumbers depending on the tableaux
are calculated for all algebras of order $n=3$ (and partially for
$n=4$). Certain classnumbers are invariants in the sense that for
algebras of same order they are equal. Algebras with several
operations of any arity are considered. A generalization of
Catalan numbers depending on homomorphisms of the structure is
proposed and corresponding generating functions set up. As an
example of this, a kind of skein polynomials are constructed
characterizing the formal build-up of iterates. As no distinction
is made between the various algebras, isomorphic or not, which
are models of the identity the results can also be formulated in
terms of varieties with signature the identity in question.
\footnote{This remark was pointed out to me by Prof. Tom C. Brown,
Simon Fraser University, Burnaby, B.C., Canada.}
\end{abstract}
\bigskip

{\bf 1. Introduction}
\par
\smallskip {\bf 1.1.} Our basic thesis is that in the ocean of
mathematical structures defined by identities the essential ones,
the ones which are liable to give rise to important developments,
are few. Such an assertion may sound presumptuous. We believe,
however, that our findings, exposed in this paper, constitute the
first steps towards validating our claim.
\par
Above statement necessitates of course to make clear the
intuitive meaning of the words structure, ocean, essential and
few. We think that every mathematician, irrespective of his
special area of interest, has the same notion in mind of the word
structure, codified in any of the various formulations which
evolved in the course of time. Ocean refers to the multitude
resulting from our freedom to imagine and write down formally any
combinations of elements and concepts entering in the definition
of structure. Structures axiomatically endowed with identities
involving operations within the underlying universe are called,
as usual {\it algebras}, resp. {\it varieties} for the totality
of algebras of same signature. The meaning of few is harder to
describe. But let us say it entails a kind of measure for a class
of algebras having a specific property in relation to the
totality of pertinent possible algebras. The establishment of
such a specific property and of the relevant measure will become
clear in the sequel at least for the algebras which will be
considered. As to the term essential structure this is actually
the great unknown and only conjectures can be made if one is to
be honest. To illustrate this characterization we would outright
say that an algebra with a binary law $Vxy$ satisfying
$VVxyz=VxVyz$ is undeniably essential. This meager functional
equation over a specified universe (set) for its variables is
simply astounding. It pervades the greater part of mathematics,
running like Ariadne's thread through the labyrinth of
mathematical theories, old and modern. The rich harvest of deep
results which have been obtained thanks to its validity - not to
say anything, platonically speaking, of those which are still out
there waiting to be unveiled - bears out the assertion. Is this
only due to the fact that its basic models are the natural
numbers and their extensions, as well as the various groups? Or
is there a much {\it deeper reason} that more complicated axiom
systems involving identities, are rare among all possible
structures we can think of? Dieudonn\'{e}, in another context,
has voiced similar doubts by saying in his introduction to the
{\it ''Abr\'{e}g\'{e} d' histoire des Mathematiques 1700-1900'' }
that ''... l' \'{e}tude des alg\'{e}bres non - associatives les
plus g\'{e}n\'{e}rales, n'ont gu\`{e}re justifi\'{e} par des
applications \`{a} des probl\`{e}mes anciens les espoirs de leurs
auteurs'' .
\par
Groups aside (associative identity), the importance of identities
for breaking new ground in mathematics can hardly be emphasized
enough. With the exception perhaps of parts of pure geometry, it
is a task indeed to seek out mathematical theories where
identities are not present either in their axioms or in their
developments. Just consider whether the theory of Lie groups and
algebras would at all have been possible were it not for the
Jacobi identity. Another case in point is the Burnside problem
for groups with generators whose elements satisfy the identity
$x^n=1$. The celebrated Euler identity in eight variables,
crucial in proving that every positive integer is a sum of four
squares, offers another example.
\par
Our investigations have been motivated by our belief,
corroborated by examination of low order algebras, that there is
a limitation restricting the number of novel algebras we can
postulate. Unless we are cunning and lucky to pick the right
ones, the ones we propose to call {\it essential}. The search for
these and their classification is the main goal. Our paper, in
four sections, is a first step in this direction. Section 2 deals
with algebras whose binary operation $Vx_{1}x_{2}$ satisfies an
identity of order $n,$ i.e. an equality between two $n$-iterates
of $V$, which cannot be {\it formally reduced} (to be explained)
to an identity of lower order. It will be shown that the number
of such algebras is of the order of $S_{n}^{2}e^{-n/16}$ as $n\to
\infty$, where $S_{n}={\frac{1}{n+1}}{2n \choose n}$, the $n$th
Catalan number. Since formally there are $S_{n}^{2} $ different
$n$-identities, this of course means that the frequency of such
algebras goes to zero as $n\to\infty$. The proof is by a novel
method we introduce namely that of tableaux, which reveal
unexpected properties yielding this result. Section 3 expounds the
semantic consequences of the defining identity. All iterates of
order higher than $n$ fall into equivalence classes of {\it
semantically equal} ones (to be explained). For the associative
identity equivalence classes of iterates of same order collapse
into a single class. Depending on the tableaux used, all
identities of order $n$ have the same classnumbers whose values
are explicitly determined as functions of certain partitions of
the Catalan numbers. If natural unicity constraints are imposed
some essential algebras have been found. For example, for $n=3$
the only essential algebras are
$VVVx_{1}x_{2}x_{3}\;x_{4}=Vx_{1}\;VVx_{2}x_{3}x_{4}$ and
$VVx_{1}Vx_{2}x_{3}\;x_{4}=Vx_{1}\;Vx_{2}Vx_{3}x_{4}$, all others
reducing to semigroups. Section 4 in form of notes and remarks,
advances the idea that the {\it method of tableaux} can be
extended to any structure provided with an equivalence relation
and generating formally new elements from elements already in the
structure. Generalizations of Catalan numbers depending on the
homomorphisms of the structure, extension of the concept of
homomorphism to several variables, skein polynomials associated
to the iterates and other allied material are included at the
end. Annexed are Exhibits of tables used in the text.
\par
\bigskip
\bigskip
{\bf 2. Formal Part}
\par
\smallskip {\bf 2.1. Ordering of iterates of a binary operation.
Construction of tableaux $A_{n}$ and $B_{n}$.} The number of
iterates of order $n$ ($n$-iterates) of a binary operation
$V(x_1x_2)$ is the same as the number of completely
parenthesizing a product of $n+1$ letters, with two factors in
each set of parentheses. As known, it is equal to the $n$th
Catalan number $S_n=\frac{1}{n+1} {2n \choose n}$. We designate
the $n$-iterates by $J_î^n(x_n\ldots x_{n+1}), i=1,\ldots ,S_n$,
putting in evidence the $n+1$ variables $x_1,\ldots ,x_{n+1}$.
\par
Suppose now we already have an ordering of the $S_{n-1}$ iterates
$Vx_{1}x_{2}$ of order $n-1$. Write them in a line in that order
and, omitting indexes, form below table by substituting
successively each variable $x$ by $Vxx$:
\begin{center}
\begin{tabular}{cccc}
$J_{1}^{n-1}(x_{1}\ldots x_{n})$ & $J_{2}^{n-1}(x_{1}\ldots
x_{n})$ & $\cdots$ & $J_{S_{n-1}}^{n-1}(x_{1}\ldots x_{n})$ \\
--------------------- & -------------------- & $\cdots$ & ----------------------- \\
$J_{1}^{n-1}(Vxx,x\ldots x)$&$J_{2}^{n-1}(Vxx,x\ldots
x)$ & $\cdots$ & $J_{S_{n-1}}^{n-1}(Vxx,x\ldots x)$ \\
$J_{1}^{n-1}(x,Vxx\ldots x)$ & $J_{2}^{n-1}(x,Vxx\ldots
x)$&$\cdots$ & $J_{S_{n-1}}^{n-1}(x,Vxx\ldots x)$ \\
$\cdots$ & $\cdots$ & $\cdots$ & $\cdots$ \\
$J_{1}^{n-1}(x,x\ldots Vxx)$ & $J_{2}^{n-1}(x,x\ldots Vxx)$&
$\cdots$ & $J_{S_{n-1}}^{n-1}(x,x\ldots Vxx)$
\end{tabular}
\end{center}

Above table comprises all $S_{n}$ $n$-iterates
$J_{i}^{n}(x_{1}\ldots x_{n+1})$, or $J_{i}^{n}$ for short,
$i=1,\ldots ,S_{n}$, each with a certain multiplicity. Tableau
$A_{n}$ (see annexed Exhibit 1, for $n=4$) is formed by assigning
to each iterate a label number from $1$ to $S_{n}$, starting with
$1$ for $J_{1}^{n-1}(Vxx,x\ldots x_{n})$ and going horizontally
through all the lines of the table, from left to right, assigning
the same label whenever the same iterate is encountered. We thus
get a well-defined ordering of the $n$-iterates, uniquely induced
by the ordering of the $(n-1)$-iterates.
\par
$n$-Iterates can be generated from $(n-1)$-iterates not just by
inserting $Vxx$ in the variable places of $(n-1)$-iterates, as
done above, but also by forming
\begin{center}
\begin{tabular}{llll}
$V(J_{1}^{n-1},x)$ & $V(J_{2}^{n-1},x)$ & $\cdots$ & $V(J_{S_{n-1}}^{n-1},x)$ \\
$V(x,J_{1}^{n-1})$ & $V(x,J_{2}^{n-1})$ & $\cdots$ &
$V(x,J_{S_{n-1}}^{n-1})$.
\end{tabular}
\end{center}

Assigning to the iterates of this table the same labels as for $A_{n}$ we
obtain tableau $B_{n}$ (see Exhibit 3).

Starting with tableaux $A_{1}$ and $B_{1}$ we can thus {\it
recursively} form tableaux $A_{n}$ and $B_{n}$, for any $n$.
Tableau $A_{n}$ has $nS_{n-1}$ entries, tableau $B_{n}$ has
$2S_{n-1}$ entries. Besides above two procedures leading to
$A_{n}$ and $B_{n}$ there is no other operation generating
$n$-iterates from $(n-1)$-iterates. There is a method based on
arithmetical properties of $A_{n}$ and $B_{n}$ to write them out,
but this is quite involved so that at this stage we omit it.
Labels have the advantage that it is much more convenient to
peruse tables with numbers and find their properties, than tables
with expressions like $VVxxVxVxx$ (see Exhibit 2). Whenever we
refer to tableaux $A_{n}$ and $B_{n}$ we shall indiscriminately
use $J_{i}^{n}$ or its label $i^{(n}$, or just $i$, if $n$ is
fixed.\smallskip

{\bf 2.2. Properties of tableaux $A_{n}$.} Let $L_{i}$,
$i=1,\ldots n$ be the $n$ sets figuring respectively in the lines
of tableau $A_{n}$. If $|M|$ denotes the cardinality of set $M$
and $M_{1}\cap M_{2}$ the intersection of $M_{1}$ and $M_{2}$ then

$$  |L_i \cap L_j|= \left\{\begin{tabular}{l@{\quad}l}
                          $S_{n-1}$ & if \quad $i=j$\\
                          0 & if \quad $|i-j|=1$\\
                          $S_{n-2}$ & otherwise
                          \end{tabular} \right.    $$

$$   |L_{i}\cap L_{j}\cap L_{k}|=\left\{\begin{tabular}{l@{\quad}l}
                       $S_{n-1}$ & if \quad $i=j=k$ \\
                       $0$ & if at least one \quad $|i-j|,\ldots ,|j-k|$\quad is $=1$ \\
                       $S_{n-3}$ & otherwise
                       \end{tabular} \right.        $$
and generally for $k=1,2,\ldots ,n$

$$   |L_{i_{1}}\cap L_{i_{2}}\cap \ldots \cap L_{i_{k}}| =\left\{
                       \begin{tabular}{l@{\quad}l}
                       $S_{n-1}$ & if $i_{1}=i_{2}=\ldots=i_{k}$ \\
                       $0$ & if at least one \quad $|i_{1}-i_{2}|,\ldots,|i_{k-1}-i_{k}|$\quad is $=1$ \\
                       $S_{n-k}$ & otherwise
                       \end{tabular} \right.      $$

$S_{n}$ are the Catalan numbers ${\frac{1}{n+1}}{2n \choose n}$.
\par
\bigskip {\bf Definition.} The {\it multiplicity} $M(J_{i}^{n})$,
$i=1,...S_{n}$ of an iterate $J_{i}^{n}$ is the number of times
it occurs in tableau $A_{n}$.
\par
The number of iterates with multiplicity $k$ is
$$T_{nk}=2^{n-2k+1}{n-1 \choose 2k-2}S_{k-1},\quad
k=1,\ldots,\left[\frac{n+1}{2}\right].$$ Above facts have been
found by induction from low values of $n$, but we have no proofs.
Maybe some competent combinatorialist can prove them.
\par
Part of the tableau $A_{n}$, starting from the right, is a Young
tableau, being a partition of $S_{n}$ into $n$ parts. The parts
are the respective lengths in each line comprising sections of
consecutive numbers of $1,2,\ldots ,S_{n}$. For example for $n=4$
(Exhibit 3) we have by $A_{4}$:

\begin{center}
\begin{tabular}{rrrrrll}
$1$ & $2$ & $3$ & $4$ & $5$ &  & $S_{n}=14=5+5+3+1$ \\
$6$ & $7$ & $8$ & $9$ & $10$ &  &  \\
&  & $11$ & $12$ & $13$ &  &  \\
&  &  &  & $14$ &  &
\end{tabular}
\end{center}

In the general case $S_{n}=c_{n,1}+c_{n,2}+\cdots +c_{n,n}$ where
the $c_{n,j}$ can be calculated recursively from following
''Pascal'' triangle

\begin{center}
\begin{tabular}{crrrrrrrcl}
$n\backslash j $ & $1$ & $2$ & $3$ & $4$ & $5$ & $\cdots$ &
& $n$ & $\cdots$    \bigskip  \\
$1 $ & $1$ &  &  &  &  &  &  &  &  $S_{1}=1$ \\
$2 $ & $1$ & $1$ &  &  &  &  &  &  &   $S_{2}=1+1=2$ \\
$3 $ & $2$ & $2$ & $1$ &  &  &  &  &  &   $S_{3}=2+2+1=5$ \\
$4 $ & $5$ & $5$ & $3$ & $1$ &  &  &  &  &   $S_{4}=5+5+3+1=14$ \\
$5 $ & $14$ & $14$ & $9$ & $4$ & $1$ &  &  &  &
$S_{5}=14+14+9+4+1=42$\\
$\cdots$ \quad & $\cdots$ & $\cdots$ & $\cdots$ & $\cdots$ &
$\cdots$ & $\cdots$ & $\cdots$ & $\cdots$   &$\cdots$ \\
$n $ & $c_{n,1}$ & $c_{n,2}$ & $\cdots$ & $\cdots$  & $\cdots$ &
$\cdots$ &
$c_{n,n}$ &  & $S_{n}=\sum_{j=1}^{n}c_{n,j}$ \\
$\cdots$ \quad & $\cdots$ & $\cdots$ & $\cdots$  &
$\cdots$ & $\cdots$ & $\cdots$ & $\cdots$ & $\cdots$ &$\cdots$  \\
\end{tabular}
\end{center}

The recursion rule is
\begin{center}
$  c_{n,j}=c_{n-1,j-1}+c_{n-1,j}+\cdots+c_{n-1,n-1}.  $
\end{center}
The numbers $c_{n,j}$ will be used in the semantic section of the
paper. An extended ''Pascal'' triangle up to $n=10$ is appended
(Exhibit 4). Out of curiosity we looked up the vertical sequences
in Sloane's Handbook of Integer Sequences and were surprised to
see that e.g. $1,3,9,28,90,297,\ldots$ are Laplace transforms
coefficients; $1,4,14,48,165,572,\ldots$ are partitions of a
polygon by a number of parts. The numbers are known in the
literature as ballot numbers, their values being
$c_{n,j}=\frac{j}{n}{2n-j-1 \choose n-1}$ (see Riordan, [1968],
Aigner [2001]).
\par
\smallskip
{\bf 2.3. Formal reducibility of an identity and incidence
matrix}. An identity $J_{i}^{n}=J_{j}^{n}$ between two
$n$-iterates is said to be of order $n$. It is {\it formally
reducible} or simply {\it reducible} to an identity of order
$n-1$ if $Vxx$ appears on both sides at the {\it same variable
place}. For example, $VVx\underline{Vxx}\;x=Vx\;V\underline{Vxx}$
$x$ of order 3 is actually formally reducible to
$VVxx\,x=Vx\,Vxx$ of order 2. Repeating the process, an identity
of order $n$ is reducible to an identity of lower order if there
is an iterate $J^{k}$, $k<n$, appearing at the same variable
place on both sides of the identity. If there is no such $J^{k}$
the identity is called {\it formally irreducible,} or simply {\it
irreducible.}
\par
In order to calculate the number of reducible resp. irreducible
identities we first form following incidence matrix relative to
tableau $A_{n}$ of all possible $S_{n}^{2}$ identities of order
$n$, including $J_{i}^{n}=J_{i}^{n}$ and counting
$J_{i}^{n}=J_{j}^{n}$ and $J_{j}^{n}=J_{i}^{n}$ as different:

\begin{center}
\begin{tabular}{llccl}
$i\backslash j$ & $J_{1}^{n}$ & $J_{2}^{n}$ & $\cdots$ &
$J_{S_{n}}^{n}$\\
\\
$J_{1}^{n}$ & $\cdots$ & $\cdots$ & $\cdots$ & $\cdots$ \\
$J_{2}^{n}$ & $\cdots$ & $\cdots$ & $\cdots$ & $\cdots$ \\
$\cdots$ & $\cdots$ & $\cdots$ & $\delta (J_{i}^{n},J_{j}^{n})$ & $\cdots$\\
$J_{S_{n}}^{n}$ & $\cdots$ & $\cdots$ & $\cdots$ & $\cdots$
\end{tabular}
\end{center}
where

$$
\delta (J_{i}^{n},J_{j}^{n})=\left\{
\begin{tabular}{l@{\quad}l}
$1$ & if \quad $J_{i}^{n}=J_{j}^{n}$ \quad reducible \\
$0$ & if \quad $J_{i}^{n}=J_{j}^{n}$ \quad irreducible.
\end{tabular} \right.   $$

The $\delta (J_{i}^{n},J_{j}^{n})$ are determined by going
successively through the lines of tableau $A_{n}$ whose
construction is such that all identities between $n$-iterates on
the same line are reducible to identities of order $n-1$ or lower.
As an example we write out the incidence matrix relative to tableau $%
A_{3}$ for all possible 25 ($S_3=5$) identities of order 3:

\begin{center}
\begin{tabular}{l@{\quad\quad}lllllr}
& $J_{1}^{3}$ & $J_{2}^{3}$ & $J_{3}^{3}$ & $J_{4}^{3}$ &
$J_{5}^{3}$ &   $\sum_{j}1$ \\
\\
$J_{1}^{3}$ & $1$ & $1$ & $0$ & $0$ & $0$ &   $2\quad$ \\
$J_{2}^{3}$ & $1$ & $1$ & $0$ & $0$ & $1$ &   $3\quad$ \\
$J_{3}^{3}$ & $0$ & $0$ & $1$ & $1$ & $0$ &   $2\quad$ \\
$J_{4}^{3}$ & $0$ & $0$ & $1$ & $1$ & $0$ &   $2\quad$\\
$J_{5}^{3}$ & $0$ & $1$ & $0$ & $0$ & $1$ &   $2\quad$\\
&  &  &  &  &  &   ---$\quad$ \\
&  &  &  &  &  &$I_{3}=11\quad$
\end{tabular}
\end{center}

The incidence matrix of order $n=4$ relative to tableau $A_{4}$
is given in Exhibit 5, and is of size $14\times 14$\quad since
$S_{4}=14.$

The sum
$$
I_{n}=\sum_{i,j=1,\ldots ,S_{n}}\delta (J_{i}^{n},J_{j}^{n})
$$

gives the total number of reducible identities and therefore
$S_{n}^{2}-I_{n}$ is the number of irreducible identities. For $n=3$,
$I_{3}=11$ as shown above. One of the main goals is to calculate $I_{n}$.
\par
\smallskip
{\bf 2.4. Connection between multiplicity of $J_{i}^{n}$  and
$\sum_{j=1}^{S_{n}}\delta (J_{i}^{n},J_{j}^{n})$.} The multiplicity
$M(J_{i}^{n})$ of $J_{i}^{n}$ was defined in subsection 2.1, as the number of
times $J_{i}^{n}$ occurs in tableau $A_{n}$. $\ \sum_{j=1}^{S_{n}}\delta
(J_{i}^{n},J_{j}^{n})$ is the sum of 1's in the $i$-th line of the incidence
matrix. We shall prove following theorem which is fundamental in our
calculation of $I_{n}$.
\par
\smallskip
{\bf Theorem.} Let $M(J_{i}^{n})=k$. Then
$$
 \sum_{j=1}^{S_{n}}\delta (J_{i}^{n},J_{j}^{n})=
{\sum_{\nu=1}^k(-1)^{\nu -1}{k \choose \nu }S_{n-\nu }}.
$$
In other words $\sum_{j=1}^{S_{n}}\delta
(J_{i}^{n},J_{j}^{n})$ does not depend on $J_{i}^{n}$ but only on the
multiplicity $k$ of $J_{i}^{n}$.
\par
\smallskip
{\bf Proof.} $J_{i}^{n}$ occurs in tableau A$_{n}$, $k$ times,
say once and only once on lines $L_{\alpha _{1}},L_{\alpha
_{2}},\ldots,L_{\alpha _{k}},$ $\ 1\le \alpha _{\nu }\le
n,~\alpha _{\nu }\neq \alpha _{\mu }$. Because of the
construction of $A_{n}$ it can never occur more than once on the
same line as all iterates of the line are formally different from
each other. Applying the inclusion - exclusion principle we get

\begin{center}
\begin{tabular}{lcl}
$\displaystyle\sum_{j=1}^{S_{n}}\delta
(J_{i}^{n},J_{i}^{n})$&$=$&$|L_{\alpha_{1}}|+|L_{\alpha_{2}}|+\cdots
+|L_{\alpha _{k}}|-$\\
 & & $|L_{\alpha _{1}}\cap L_{\alpha_{2}}| -\cdots -
 |L_{\alpha _{k-1}}\cap L_{\alpha _{k}}| +$ \\  &  &  \\
& &$|L_{\alpha_{1}}\cap L_{\alpha _{2}}\cap L_{\alpha _{3}}|
+\cdots+ |L_{\alpha_{k-2}}\cap L_{\alpha _{k-1}}\cap L_{\alpha
_{k}}|-$\\  &  &  \\
&  & $\cdots +(-1)^{k-1}|L_{\alpha_{1}}\cap L_{\alpha _{2}}\cap
\cdots \cap L_{\alpha _{k}}|.$
\end{tabular}
\end{center}
\par
The case of one of the $|L_{\alpha _{\nu }}\cap L_{\alpha
_{\mu}}|$ being $0$ is impossible. This would mean according to
2.2 that $J_{i}^{n}$ occurs in tableau $A_{n}$ in two consecutive
lines. It can be proved however by induction from $A_{n-1}$ to
$A_{n}$ that this cannot happen. In fact $A_{n}$ is a
''Schachtelung'' of Young tableaux distanced vertically by at
least 2 lines from each other as seen e.g. for $A_{n}$ on annexed
Exhibit 3. As a consequence all $\left| L_{\alpha _{1}}\cap
L_{\alpha _{2}}\cap \ldots \right| $ are also different from $0$.
Substituting their values as per 2.2 considering that $\alpha
_{\nu }\neq \alpha _{\mu }$ we obtain $\sum_{j=1}^{S_{n}}\delta
(J_{i}^{n},J_{j}^{n})=\sum_{\nu =1}^{k}(-1)^{\nu -1}{k \choose
\nu }S_{n-\nu }$ as required.
\par
\smallskip
{\bf 2.5. Calculation of $I_{n}$.} We can now calculate
$I_{n}=\sum_{i=1}^{S_{n}}\sum_{j=1}^{S_{n}}\delta (J_{i}^{n},J_{j}^{n})$.
Since $\sum_{j=1}^{S_{n}}\delta (J_{i}^{n},J_{j}^{n})$ is the number of
times $J_{i}=J_{j}$ is reducible in the $i$-th line of the incidence matrix
and, as indicated in 2.2, there are
$$
T_{n,k}=2^{n-2k+1}{n-1 \choose 2k-2}S_{k-1}\ \ \ \ \ \
{n-1\choose2k-2} \ge 1
$$
iterates $J_{i}$ with multiplicity $k$ we get for $I_{n}$:
$$I_{n}=\sum_{k=1}^{[\frac{n+1}{2}] }T_{n,k}\left (
\sum_{j=1}^{S_{n}}\delta (J_{i}^{n},J_{j}^{n})\right).$$ The upper
limit $[\frac{n+1}{2}] $ for $k$ is obtained from $n-1\geq 2k-2$,
all $T_{n,k}$ with $k>[\frac{n+1}{2}] $ being zero.
\par
Inserting for $\sum_{j=1}^{S_{n}}\delta (J_{i}^{n},J_{j}^{n})$ its value
found in 2.4 we get
$$
I_{n}=\sum_{k=1}^{[\frac{n+1}{2}]}T_{n,k}\left[{k \choose 1}
S_{n-1}-{k \choose 2}S_{n-2}+\cdots +(-1)^{k-1}{k \choose
k}S_{n-k}\right]
$$
and by changing the order of summation
$$I_{n}=\sum_{k=1}^{\frac{n+1}{2}}(-1)^{k-1}S_{n-k}\left[
{k \choose k}T_{n,k}+{k+1 \choose k}T_{n,(k+1)}+\cdots +
{\left[\frac{n+1}{2}\right] \choose k}T_{n,
\left[\frac{n+1}{2}\right]}\right].
$$
\par
 The evaluation of the brackets in this formula has transcended our efforts
but straightforward calculations of $I_{n}$ based on tableaux $A_{n}$ for
$n=3,4,5,6,7$ and their respective incidence matrices have led us to surmise
that
$$\sum_{\nu =0}^{[\frac{n+1}{2}]-k}{k+\nu \choose k}T_{n,k+\nu }
= {n-k+1 \choose k}S_{n-k}.$$
For $k=1$ this formula reduces to

$$
1T_{n,1}+2T_{n,2}+\cdots+\left[
\frac{n+1}{2}\right]T_{n,\left[\frac{n+1}{2}\right]}= nS_{n-1},
$$
which is certainly correct considering that tableau $A_{n}$ has
$nS_{n-1}$ entries. Maybe application of the ''Snake Oil Method''
or the ''Wilf-Zeilberger Method'' which we have not tried will
succeed in proving this identity (see Wilf [1990]).
\par
Substituting the brackets in above formula by their values
${n-k+1 \choose k}S_{n-k}$ and skipping the upper limit $[{\frac{n+1}{2}}]$ in the
summation sign, as all ${n-k+1 \choose k}$ are zero for
$k>[ {\frac{n+1}{2}}]$, the expression for $I_{n}$ takes finally the form
$$
I_{n}=\sum_{k=1}^{\infty }(-1)^{k-1}{n-k+1 \choose k}S_{n-k}^{2}.
$$
\par
\smallskip
{\bf 2.6. Asymptotic evaluation of $I_{n}$.} As there are $S_{n}^{2}$ formally different algebras, counting them as in section 2.3,
defined by an identity $J_{i}^{n}=J_{j}^{n}$, and the number of reducible
algebras is $I_{n}$,
$$
\frac{I_{n}}{S_{n}^{2}}=\sum_{k=1}^{\infty }(-1)^{k-1}{n-k+1
\choose k }\left( \frac{S_{n-k}}{S_{n}}\right) ^{2} $$ is the
probability (measure) for an algebra to be reducible and hence
$$
1-\frac{I_{n}}{S_{n}^{2}} $$

is the probability that it be irreducible.
\par
From the recurrence $S_{n}=2\frac{2n-1}{n+1}S_{n-1}$ for the
Catalan numbers it can be easily shown by induction from $n$ to
$n+1$ that
\begin{center}
$\frac{S_{n-k}}{S_{n}}\to \frac{1}{4^{k}} \quad$ for \quad  $n\to
\infty .$
\end{center}

\par
On the other hand ${{n-k+1} \choose {k}}$ is $\sim
\frac{n^{k}}{k!} \quad$ for $ n\to \infty $. The general term
therefore of above series behaves like

$$ (-1)^{k-1}\frac{1}{k!}\frac{1}{4^{2k}}n^{k}. $$

Summing over $k$ we can approximate $I_{n}/S_{n}^{2}$ by

$$
\sum_{k=1}^{\infty }(-1)^{k-1}\frac{1}{k!}\left( \frac{n}{16}\right)
^{k}=1-e^{-n/16},
$$
and hence,
$$\lim_{n\to \infty }\frac{I_{n}}{S_{n}^{2}}\sim
\lim_{n\to \infty }\left( 1-e^{-n/16}\right) =1. $$
\par
Above reasoning is of course not rigorous needing either a direct
evaluation of the series or a "Fehlerabsch{\"a}tzung" of the
difference with $1-e^{-n/16}$.
\par
\smallskip
Stronger results can be obtained if instead of counting only the
reducible identities of tableau $A_{n}$ we do also take into
account the reducible identities of tableau $B_{n}$, which do not
already occur in tableau $A_{n}$ . For example from tableau
$B_{4}$ (Exhibit 3) we infer that $ J_{1}^{4}=J_{8}^{4}$ and
$J_{9}^{4}=J_{13}^{4}$ are formally reducible, a fact {\it which
can not be deduced from tableau} $A_{4}$. By ''adding'' tableau
$A_{n}$ and $B_{n}$ to form $ A_{n}\oplus B_{n}$ with
$(n+2)S_{n-1}$ entries, we can again evaluate its corresponding
incidence matrix which remarkably has the same properties as the
incidence matrix of tableau $A_{n}$. The theorem that the sum of
1's on a line of the incidence matrix is the same for all
$J_{i}^{n}$ having the same multiplicity remains unchanged. The
occurrence, however, of an iterate with multiplicity $k$ has now
been found to be
$$
T_{n,k}^{A_{n}\oplus B_{n}}=T_{n,k}+2\left(
T_{n-1,k-1}-T_{n-1,k}\right) ,
$$
with obvious notation. We have not pursued the calculations to find an
expression for
$$
I_{n}^{A_{n}\oplus B_{n}}=\sum_{k=1}^{\left[\frac{n+1}{2}\right]
}T_{n,k}^{A_{n}\oplus B_{n}} \left[ {k \choose 1}S_{n-1}-{k
\choose 2}S_{n-2}+\cdots +(-1)^{k-1} {k \choose k}S_{n-k}\right],
$$ but evaluations up to $n=7$ show as expected a much faster
convergence of $\ I_{n}^{A_{n}\oplus B_{n}}/S_{n}^{2}$ to $1$,
than that of $I_{n}/S_{n}^{2}$.
\par
\smallskip
{\bf 2.7. Excursus on the implications of
$I_{n}/S_{n}^{2}\rightarrow 1$.} The fact that
$I_{n}^{A_{n}}/S_{n}^{2}\rightarrow 1$ shows that almost all
$n$-identities reduce formally to identities of lower order. Since
$$  \frac{S_{n}^{2}-I_{n}}{S_{n}^{2}}\sim e^{-n/16}$$

this means that {\it there is almost no formally irreducible identity as}
$n\rightarrow \infty $. This explicates the word {\it few} we used in the
introduction regarding the possibility of stipulating structures
axiomatically provided with an identity.
\par
{\it Formal irreducibility is an imperative requirement} of an
identity figuring in an axiom system: nobody would define a
semigroup by $ VVVx_{1}x_{2},x_{3},x_{4}=VVx_{1}x_{2},Vx_{3}x_{4}$
($Vx_1x_2$ occurring in the same variable space on both sides) and
not by $ VVx_{1}x_{2},x_{3}=Vx_{1},Vx_{2}x_{3}$ . If we ask
therefore a mathematician to postulate an algebra with an
identity of say order $n=35$ and claim it is a new structure,
chances are high that his new product is redundant, as being
reducible to a lower order, maybe even to the associative law. If
the question to turn out a kilometer long identity is put to a
computer, theoretically speaking, it is almost certain that the
brainchild will be a hydrocephalus.
\par
At this point we cannot resist the temptation to digress somehow from the
mainstream of our exposition. We have always wondered why mathematical axiom
systems, theorems and formulas are so short, rarely taking more than a
quarter page of a book to be {\it formally} expressed within a certain
alphabet. Obviously, this is due to our anthropomorphic yardstick which is
limiting our capacity to consider and manipulate long expressions. With the
help of computers the limits are just receding farther off. Unavoidably,
however, this raises a major question which, to our knowledge we haven't
seen addressed to before. Are there mathematical ''truths'' which would
require at least, say $100^{100^{100}}$ pages of an ordinary book to be
expressed in some formal language, not including proofs? Clearly,
abbreviations like $\sum $ and formulations depending on a parameter $n$ and
increasing in length with increasing $n$ cannot be considered to yield such
statements if a specific large number is inserted for $n$. And of what
nature and content would such monster theorems be which we wouldn't be able
to read, let alone grasp?
\par
But back to terra ferma and our smallish theorems until such a
monster is discovered - or not.\bigskip

{\bf 3. Semantic Part}
\par
\smallskip
{\bf 3.1. Classes of semantically equal iterates resulting from
an identity.} Semantic equality means that in an identity
$J_i^n(x_1\ldots x_{n+1})=J_j^n(x_1\ldots x_{n+1})$ the two sides
are formally different but are the same (equal) for all $x_\nu $
running over the prescribed space for the variables. The usual
properties of an equivalence relation hold.
\par
The results of section 2 were mainly obtained by investigating
the lines of tablaux $A_{n}$ (occasionally also those of tablaux
$B_{n}$). In the following, we will investigate the columns of
tableaux $A_{n}$ and $B_{n}$ and derive semantic information
regarding the iterates of $Vx_{1}x_{2}$ of order higher than $n$,
as a consequence of the defining identity.
\par
To illustrate the general case we will take an example drawn from
tableaux $A_{4}$ and $B_{4}$ (Exhibit 3). Suppose we postulate a
3-algebra defined by the identity $J_{2}^{3}=J_{4}^{3}$ or
written out $V(Vxx,\,Vxx)=V(x,\,V(Vxx,\,x)).$ Using the
enumerative labeling described in 2.1,
$(2^{(3},4^{(3},10^{(4},\ldots$ stand respectively for
$J_{2}^{3},J_{4}^{3},J_{10}^{4},\ldots )$ and going down columns 2
and 4 of tableaux $A_{4}$ and $B_{4}$ we deduce following
identities:

\begin{center}
\begin{tabular}{ccc}
From Tableau $A_{4}$ &  & From Tableau $B_{4}$ \\
$2^{(3}=4^{(3}$ &  & $2^{(3}=4^{(3}$ \\
--------------------- &  & --------------------- \\
$2^{(4}=4^{(4}$ &  & $3^{(4}=8^{(4}$ \\
$7^{(4}=9^{(4}$ &  & $10^{(4}=13^{(4}$ \\
$4^{(4}=12^{(4}$ &  &  \\
$5^{(4}=10^{(4}$ &  &
\end{tabular}
\end{center}
\par
This means that the 14 iterates of order 4 fall into 10 classes of
semantically equal iterates if tableau $A_{4}$ is used, namely (skipping
upper indices):

\begin{center}
\{1\}, \{2, 4, 12\}, \{3\}, \{5, 10\}, \{6\}, \{7, 9\}, \{8\}, \{11\},
\{13\}, \{14\}.
\end{center}
\par
If $B_{4}$ is also taken into account the number of classes is
reduced to 8 because of the equivalence properties of semantic
equality:

\begin{center}
\{1\}, \{2, 4, 12\}, \{3, 8\}, \{5, 10, 13\}, \{6\}, \{7, 9\}, \{11\},
\{14\}.
\end{center}

We can now repeat the process by using either of above sets as
defining identities and obtain via tableaux $A_5$ and $B_5$ the
equivalence classes of order 5 and their classnumbers (number of
classes). This procedure, which can be continued indefinitely,
\textit{splits the set of all iterates of any order into
equivalence classes of semantically equal ones,} as a consequence
of the initial defining identity $2^{(3}=4^{(3}$.
\par
In the general case we start with an algebra postulated by an
identity $J^{n}=\overline{J^{n}}$ between two formally different
iterates - several identities can also be considered as we shall
see in section 4. We can apply above procedure by using tableaux
$A_{n+k}$ and $B_{n+k}$, $k\geq 1$, in a certain prescribed
order, in which case the corresponding class sets and
classnumbers will depend on the chosen order. We shall restrict
ourselves to using exclusively either tableaux $A_{n+k}$ or both
tableaux $A_{n+k}$ and $B_{n+k}$, in conjunction, which as said in
2.6 we denote by $A_{n+k}\oplus B_{n+k}$. In the first case we
have found a formula to calculate the classnumbers for any $k$,
(see following theorem), without actually going through the
laborious procedure of forming the successive class sets. This
however is not giving the best results as shown by above example.
The actual minima of the classnumbers are obtained by using
consistently tableaux $A_{n+k}\oplus B_{n+k}$. We shall denote by
$ h_{n+k}^{A},h_{n+k}^{B}$ and $h_{n+k}^{A\oplus B}$ respectively
the classnumbers in reference to the tableaux employed. As
$h_{n+k}^{A\oplus B}$ is the minimal number of the equivalence
classes we shall omit the indication $A\oplus B$ and simply call
$h_{n+k}=h_{n+k}(J^{n}=\overline{J^{n}})$ the classnumbers of the
identity, for $k \ge 0$. Below $n$ the classnumbers are of course
the Catalan numbers since the identity is of order $n$.

Writing $H_{i}^{n+k}$, $i=1,...,h_{n+k}$ for the classes of order
$n+k$ we can list them as follows together with their $h_{n+k}$:

\begin{center}
\begin{tabular}{lllllllllll}
$H_{1}^{0}$ &  &  &  &  &  &  &  &  &  & $h_{0}=1$ \\
$H_{1}^{1}$ &  &  &  &  &  &  &  &  &  & $h_{1}=1$ \\
$H_{1}^{2}$ & $H_{2}^{2}$ &  &  &  &  &  &  &  &  & $h_{2}=2$ \\
$H_{1}^{3}$ & $H_{2}^{3}$ & $H_{3}^{3}$ & $H_{4}^{3}$ & $H_{5}^{3}$ &  &  &
&  &  & $h_{3}=5$ \\
$\cdots$ & $\cdots$ & $\cdots$ & $\cdots$ & $\cdots$ &  &  &  &
&  &  \quad $\cdots$\\
$H_{1}^{n-1}$ & $H_{2}^{n-1}$ & $\cdots$ & $\cdots$ & $\cdots$ &
$H_{S_{n-1}}^{n-1}$ &  &
&  &  & $h_{n-1}=S_{n-1}$ \\
$H_{1}^{n}$ & $H_{2}^{n}$ & $\cdots$ & $\cdots$ & $\cdots$ &
$\cdots$ & $H_{S_{n}-1}^{n}$ &  &
&  & $h_{n}=S_{n}-1$ \\
$\cdots$ & $\cdots$ & $\cdots$ & $\cdots$ & $\cdots$ & $\cdots$ &
$\cdots$ & $\cdots$ & $\cdots$ &
$\cdots$ & \quad $\cdots$ \\
$H_{1}^{n+k}$ & $H_{2}^{n+k}$ & $\cdots$ & $\cdots$ & $\cdots$ &
$\cdots$ & $\cdots$ & $\cdots$ & $H_{h_{n+k}}^{n+k}$ &  &
$h_{n+k}$ \\
$\cdots$ & $\cdots$ & $\cdots$ & $\cdots$ & $\cdots$ & $\cdots$ &
$\cdots$ & $\cdots$ & $\cdots$ & $\cdots$ & \quad $\cdots$
\end{tabular}
\end{center}

\par
The main objective is to calculate the classnumbers $h_{n+k}$ for
a given identity of order $n$ and to classify the identities
according to their classnumbers. A major find, conjectured,
towards this aim, if we restrict ourselves to using only the
sequence of the tableaux $A_{1},A_{2},\ldots,A_{k},\ldots $ is the
\par
\smallskip
{\bf Theorem.} All identities $J^{n}=\overline{J^{n}}$ of order
$n$ have the same classnumber $h_{n+k}^{A}$ namely,
$$
h_{n+k}^{A}(J^{n}=\overline{J^{n}})=S_{n+k}-c_{n+k+1,n+1}\qquad k\geq 0
$$
where the $c_{i,,j}$ are the numbers of the ''Pascal'' triangle
defined in 2.2.
\par
\bigskip Since $h_{n+k}<h_{n+k}^{A}$, as any equivalence class from
tableau $A\oplus B $ comprises all iterates of the corresponding
equivalence class of tableau $A $, we have
$$
h_{n+k}(J^{n}=\overline{J^{n}})<S_{n+k}-c_{n+k+1,n+1}\qquad k\geq 0
$$
which of course is stronger than the trivial $h_{n+k}<S_{n+k}$,
especially for big $k$.
\par
Using tableaux $B_{n}$ only it is easily seen that $
h_{n+k}^{B}=S_{n+k}-2^{k} $, as all entries in the 2 lines of
$B_{n+k}$ appear only once.
\par
\smallskip
{\bf 3.2. The algebra of classes.} The classes $H_{i}^{p}$ can be
taken as elements of a new algebra with a binary composition law
$W(H_{i}^{p},H_{j}^{q})$ defined as follows:
\par
\smallskip
{\bf Definition.} If $J_{i}^{p}$ is a representative of the class
$H_{i}^{p}$ and $J_{j}^{q}$ is a representative of the class
$H_{j}^{q}$ then the class $W(H_{i}^{p},H_{j}^{q})$ is defined as
the class of order $p+q+1$ to which the iterate
$V(J_{i}^{p},J_{j}^{q})$ of order $p+q+1$ belongs.
\par
The class $W(H_{i}^{p},H_{j}^{q})=H_{\lambda }^{p+q+1}$ is
uniquely defined as any other pair $J_{\mu }^{p},~J_{\nu }^{q}$,
of representatives will give an iterate $V(J_{\mu }^{p},J_{\nu
}^{q})$ semantically equal to $V(J_{i}^{p},J_{j}^{q})$. The new
law $W(H_{1},H_{2})$ defined over the classes satisfies the same
identity as the original algebra. It is a L\"{o}wenheim - Skolem
denumerable model of the identity. The classnumbers $h_{n+k}$ now
play the same role as the Catalan numbers for the original
algebra and corresponding tableaux can be formed. The process can
be repeated giving a series of algebras.
\par
\smallskip
{\bf 3.3 Semantically reducible algebras.} By {\it semantic
reducibility} we mean that $V(x,y)$ besides satisfying the
defining formally irreducible identity, does also satisfy,
because of inherent reasons of the structure, an identity of lower
order. In the opposite case the identity is called {\it
semantically irreducible} or {\it essential.}
\par
So far the algebras considered were subjected to no constraint
whatsoever besides the axiom of the defining identity. In the
rest of this subsection we shall additionally provide that if
$V(x,a)=V(y,a)$ then $x=y$ and similarly
$V(a,x)=V(a,y)\rightarrow x=y$. These are natural restrictions
and practically fulfilled in many structures of importance
(unicity of solutions of $V(x,a)=b$ resp. $V(a,x)=b$ in case they
exist).
\par
We have made a thorough investigation of all 3-algebras subject to above
unicity property of $V\,xy$ which led to the determination of all essential
structures of that order as well as of their classnumbers.
\par
The following theorem summarizes these results and the proof
straightforward but tedious, clarifies the concept of semantic
reducibility.
\par
\smallskip
{\bf Theorem.} There are only 2 algebras of order 3, semantically
irreducible i.e.
$$
V(V(Vx_{1}x_{2},x_{3}),x_{4})=V(x_{1},V(Vx_{2}x_{3},x_{4})),\qquad
(1^{(3}=4^{(3})  $$ and
$$
V(V(x_{1},Vx_{2}x_{3}),x_{4})=V(x_{1},V(x_{2},Vx_{3}x_{4})),\qquad
(3^{(3}=5^{(3}).
$$
\par
All other 3-algebras are reducible, formally or semantically to
algebras of lower order, i.e. to semigroups or to the free
algebra of $Vxy$ with no identity. This theorem is true provided
$V(x,y)$ fulfills above unicity conditions.
\par
\smallskip
{\bf Proof.} We form the incidence matrix (see 2.3) $5\times 5$ of
all 3-order identities resulting {\it from both tableaux }$A_{3}$
{\it and} $B_{3}$, i.e. from $A_{3}\oplus B_{3}$, in order to
include also the reducible identities which result from tableau
$B_{3}$ but not from tableau $A_{3}$:

\begin{center}
\begin{tabular}{llllll}
& $J_{1}$ & $J_{2}$ & $J_{3}$ & $J_{4}$ & $J_{5}$ \\
$J_{1}$ & $1$ & $1$ & $1$ & $0$ & $0$ \\
$J_{2}$ &  & $1$ & $0$ & $0$ & $1$ \\
$J_{3}$ &  &  & $1$ & $1$ & $0$ \\
$J_{4}$ &  &  &  & $1$ & $1$ \\
$J_{5}$ &  &  &  &  & $1$%
\end{tabular}
\end{center}

Because of symmetry only the upper triangle is considered.
\par
We see that the only identities which are not formally reducible
are the following five, indicated by $0$'s in the matrix:
$$
1^{(3}=4^{(3},\quad 1^{(3}=5^{(3},\quad 2^{(3}=3^{(3},\quad
2^{(3}=4^{(3},\quad 3^{(3}=5^{(3}.
$$
\par
Consider f.e. $1^{(3}=5^{(3}$ and construct the equivalence
classes of iterates of order 4 and 5, by applying successively
tableaux $A_{4}\oplus B_{4}$ and $A_{5}\oplus B_{5}$ as outlined
in 3.1. From tableau $A_{4}\oplus B_{4}$ we get successively
(writing $i^{(n}$ or $i$ for $J_{i}^{n}$), by collection of
classes having common elements, following 6 identities which
yield 8 $(h_{4}=8)$ equivalence classes for the 14 $(S_{4}=14)$
iterates of order 4:

\begin{center}
\begin{tabular}{l}
$1^{(3}=5^{(3}$ \\
------------ \\
$1^{(4}=5^{(4}$ \\
$6^{(4}=10^{(4}$ \\
$3^{(4}=13^{(4}$ \\
$2^{(4}=14^{(4}$ \\
$1^{(4}=11^{(4}$ \\
$9^{(4}=14^{(4}$
\end{tabular}
\end{center}
\begin{center}
\{1,5,11\},\quad (2,9,14\},\quad \{3,13\},\quad \{4\},\quad
\{6,10\},\quad \{7\},\quad \{8\},\quad \{12\}.
\end{center}
\par
Repeating the same procedure with above classes we get from
tableaux $A_{5}\oplus B_{5}$ (Exhibit 3) -singletons can be
omitted- following 14 classes for the 5-iterates (the equality
sign between the resulting 5-iterates has been omitted)

\begin{center}
\begin{tabular}{ccccccc}
$1^{(4}=5^{(4}=11^{(4}$ &  & $2^{(4}=9^{(4}=14^{(4}$ &  & $3^{(4}=13^{(4}$ &
& $6^{(4}=10^{(4}$ \\
--------------------- &  & ----------------------- &  & -------------- &  &
-------------- \\
$~~1\qquad ~~5\qquad 11$ & \multicolumn{1}{l}{} & $~~2\qquad ~~9\qquad 14$ &
\multicolumn{1}{l}{} & $~~3\qquad 13$ & \multicolumn{1}{l}{} & $~~6\qquad 10$
\\
$15\qquad 19\qquad 25$ & \multicolumn{1}{l}{} & $16\qquad 23\qquad 28$ &
\multicolumn{1}{l}{} & $17\qquad 27$ & \multicolumn{1}{l}{} & $20\qquad 24$
\\
$~~6\qquad 10\qquad 34$ & \multicolumn{1}{l}{} & $~~7\qquad 32\qquad 37$ &
\multicolumn{1}{l}{} & $~~8\qquad 36$ & \multicolumn{1}{l}{} & $29\qquad 33$
\\
$~~3\qquad 13\qquad 38$ & \multicolumn{1}{l}{} & $~~4\qquad 26\qquad 41$ &
\multicolumn{1}{l}{} & $11\qquad 40$ & \multicolumn{1}{l}{} & $17\qquad 27$
\\
$~~2\qquad 14\qquad 30$ & \multicolumn{1}{l}{} & $~~5\qquad 24\qquad 42$ &
\multicolumn{1}{l}{} & $~~7\qquad 37$ & \multicolumn{1}{l}{} & $16\qquad 28$
\\
$~~1\qquad 11\qquad 29$ & \multicolumn{1}{l}{} & $~~3\qquad 22\qquad 38$ &
\multicolumn{1}{l}{} & $~~6\qquad 34$ & \multicolumn{1}{l}{} & $15\qquad 25$
\\
$23\qquad 28\qquad 39$ & \multicolumn{1}{l}{} & $24\qquad 36\qquad 42$ &
\multicolumn{1}{l}{} & $26\qquad 41$ & \multicolumn{1}{l}{} & $32\qquad 37$%
\end{tabular}
\end{center}

\begin{center}
\{1, 5, 8, 11, 20, 24, 29, 33, 36, 40, 42\}, \{2, 9, 14, 30\},
\{3, 13, 22, 38\}

\{4, 26, 41\}, \ (6, 10, 34\}, \{7, 32, 37\}, \{12\}, \{15, 19,
25\}, \{16, 23, 28, 39\}

\{17, 27\}, \{18\}, \{21\}, \{31\}, \{35\}
\end{center}
\par
As we see one of the 14 $5$-order iterate classes is (1, 5, 8, 11,
20, 24, 29, 33, 36, 40, 42) which contains both $8^{(5}$and
$11^{(5}$ as semantically equal. Let us write out in full the
identity $8^{(5}=11^{(5}$ by substituting for the labels their
corresponding 5-iterates (see Exhibit 2):
$$
V(V(Vx_{1}x_{2},\,V(Vx_{3}x_{4},\,x_{5})),\,x_{6})=V(V(Vx_{1}x_{2},\,V(x_{3},
\,Vx_{4}x_{5})),\,x_{6}).
$$
Because of the unicity assumption stipulated in the beginning of
this section, we get successively
$$
V(Vx_{1}x_{2},\,V(Vx_{3}x_{4},\,x_{5}))=V(Vx_{1}x_{2},\,V(x_{3},
\,Vx_{4}x_{5}))
$$
$$ V(Vx_{3}x_{4},\,x_{5})=\,V(x_{3},\,Vx_{4}x_{5}).$$
\par
The last identity is of course the associative law so that the
algebra $1^{(3}=5^{(3}$ which is formally irreducible, is
actually a semigroup.
\par
The same can be shown to be true for the algebras $2^{(3}=3^{(3}$
and $2^{(3}=4^{(3}$, leaving as only semantically irreducible
(essential!) algebras, the cases $1^{(3}=4^{(3}$,
$3^{(3}=5^{(3}$, as stated in the theorem. We would like to
stress that above theorem is true at the present stage of the
investigations as for the proof it suffices to form the
equivalence classes of order up to 5. But who knows whether by
proceeding to equivalence classes of orders $6$, $7$, $\ldots$,
$n_{0}$, we won't find a pair of $n_{0}$ -iterates which would
similarly also eliminate the cases $1^{(3}=4^{(3}$, and
$3^{(3}=5^{(3}$? We have arithmetic evidence that this is rather
unlikely.
\par
It is possible to spare to trouble of writing out the iterates,
which is extremely cumbersome, if we know beforehand which pairs
of $n$-iterates automatically imply the semantic equality of
iterates of order $2,3,...,n-1$, respectively. This can be done
with the help again of tableau $A_{n}$. The number of such pairs
can be shown to be
$$
(n-k+1)S_{n-k}\frac{S_{k}(S_{k}-1)}{2}\qquad 2\leq k\leq n-1.
$$
\par
Annexed here-with (Exhibit 6), by way of examples, are the tables
of such iterates for $n=5$, $k=2$ (couples) and for $n=5$, $k=3$
(triples). They figure respectively in the two tableaux $A_5$,
connected by the sign $\updownarrow $ and should be read off
vertically. Thus in tableaux $A_5$, $8^{(5}$ and $11^{(5}$ figure
in same column connected by $\updownarrow $, implying
associativity, as shown above. For $n=5$, $k=3$ all pairs of
$5$-order iterates in same column connected by $\updownarrow $, if
semantically equal, automatically implicate the semantic equality
of two 3-order iterates.
\par
The rules of placement of the relational sign $\updownarrow $, which we
won't expound, involve again the Catalan numbers and the numbers $c_{nj}$ of
the ''Pascal'' triangle defined in 2.2. This is further evidence of the
wealth of information encoded in the tableaux $A_n$.
\par
\bigskip {\bf 3.4. The series of classnumbers $h_{n}$ for
identities of order 3 and 4.} Equality in this subsection means
semantic equality.
\par
As shown in 3.1 the partition of all $N$-iterates into $h_{N}$
classes of equal iterates induces as a result of tableaux
$A_{N+1}\oplus B_{N+1}$ a partition of all $(N+1)$-iterates into
$h_{N+1}$classes of equal iterates. The problem of one defining
$n$-identity, is starting with $h_{n}=S_{n}-1,$ to calculate the
series $h_{n+1},...,h_{n+k},...$ . Generating functions $
\sum_{0}^{\infty }h_{n}t^{n}$ can be defined and calculated in
some cases, for instance for $3^{(3}=5^{(3}$.
\par
We have carried out the investigation for all 3-identities and
calculated the respective classnumbers up to $n=7$, as well as
the number of classes containing only 1 element (singletons).
Beyond that the process gets unwieldy because of the fast growth
of the Catalan numbers $S_{n}$. The following list summarizes our
results for the 10 possible identities of order 3, obviously
excluding $J_{i}^{(3}=J_{i}^{(3}$ and counting $
J_{i}^{(3}=J_{j}^{(3}$ and $J_{j}^{(3}=J_{i}^{(3}$ as one
identity. We have also included the results for the formally
reducible algebras because of the surprising fact that for some
of them the number of singletons are the Fibonacci numbers
$F_{n}$:
\newpage
\begin{center}
\begin{tabular}{ccccccc}
Algebra &  & Formally &  & Classnumber &  & Number of \\
&  & reducible/ &  & $h_{n}$ &  & singletons \\
&  & irreducible &  &  &  &  \\
\multicolumn{1}{l}{-----------} & \multicolumn{1}{l}{} & \multicolumn{1}{l}{
--------------} & \multicolumn{1}{l}{} & \multicolumn{1}{l}{
--------------------} & \multicolumn{1}{l}{} & \multicolumn{1}{l}{
--------------} \\
\multicolumn{1}{l}{$1^{(3}=2^{(3}$} & \multicolumn{1}{l}{} &
\multicolumn{1}{l}{reducible} & \multicolumn{1}{l}{} & \multicolumn{1}{l}{$%
2^{n-1}\qquad n\geq 3$} & \multicolumn{1}{l}{} & \multicolumn{1}{l}{$%
F_{n}\quad n\geq 3$} \\
\multicolumn{1}{l}{$1^{(3}=3^{(3}$} & \multicolumn{1}{l}{} &
\multicolumn{1}{l}{reducible} & \multicolumn{1}{l}{} & \multicolumn{1}{l}{$%
2^{n-1}\qquad n\geq 3$} & \multicolumn{1}{l}{} & \multicolumn{1}{l}{$%
F_{n}\quad n\geq 3$} \\
\multicolumn{1}{l}{$1^{(3}=4^{(3}$} & \multicolumn{1}{l}{} &
\multicolumn{1}{l}{irreducible} & \multicolumn{1}{l}{} & \multicolumn{1}{l}{$%
2^{n-1}\qquad n\geq 3$} & \multicolumn{1}{l}{} & \multicolumn{1}{l}{$n\quad
~~\,n\geq 3$} \\
\multicolumn{1}{l}{$1^{(3}=5^{(3}$} & \multicolumn{1}{l}{} &
\multicolumn{1}{l}{irreducible} & \multicolumn{1}{l}{} & \multicolumn{1}{l}{$%
4\qquad \quad ~\,n=3$} & \multicolumn{1}{l}{} & \multicolumn{1}{l}{$n\quad
~~\,n\geq 3$} \\
\multicolumn{1}{l}{} & \multicolumn{1}{l}{} & \multicolumn{1}{l}{} &
\multicolumn{1}{l}{} & \multicolumn{1}{l}{$8\qquad \quad ~\,n=4$} &
\multicolumn{1}{l}{} & \multicolumn{1}{l}{} \\
\multicolumn{1}{l}{} & \multicolumn{1}{l}{} & \multicolumn{1}{l}{} &
\multicolumn{1}{l}{} & \multicolumn{1}{l}{$14\qquad ~~~n=5$} &
\multicolumn{1}{l}{} & \multicolumn{1}{l}{} \\
\multicolumn{1}{l}{} & \multicolumn{1}{l}{} & \multicolumn{1}{l}{} &
\multicolumn{1}{l}{} & \multicolumn{1}{l}{$20\qquad ~~~n=6$} &
\multicolumn{1}{l}{} & \multicolumn{1}{l}{} \\
\multicolumn{1}{l}{} & \multicolumn{1}{l}{} & \multicolumn{1}{l}{} &
\multicolumn{1}{l}{} & \multicolumn{1}{l}{$16\qquad ~~~n=7$} &
\multicolumn{1}{l}{} & \multicolumn{1}{l}{} \\
\multicolumn{1}{l}{$2^{(3}=3^{(3}$} & \multicolumn{1}{l}{} &
\multicolumn{1}{l}{irreducible} & \multicolumn{1}{l}{} & \multicolumn{1}{l}{$%
4\qquad \quad ~\,n=3$} & \multicolumn{1}{l}{} & \multicolumn{1}{l}{$n\quad
~~\,n\geq 3$} \\
\multicolumn{1}{l}{} & \multicolumn{1}{l}{} & \multicolumn{1}{l}{} &
\multicolumn{1}{l}{} & \multicolumn{1}{l}{$8\qquad \quad ~\,n=4$} &
\multicolumn{1}{l}{} & \multicolumn{1}{l}{} \\
\multicolumn{1}{l}{} & \multicolumn{1}{l}{} & \multicolumn{1}{l}{} &
\multicolumn{1}{l}{} & \multicolumn{1}{l}{$14\qquad ~~~n=5$} &
\multicolumn{1}{l}{} & \multicolumn{1}{l}{} \\
\multicolumn{1}{l}{} & \multicolumn{1}{l}{} & \multicolumn{1}{l}{} &
\multicolumn{1}{l}{} & \multicolumn{1}{l}{$20\qquad ~~~n=6$} &
\multicolumn{1}{l}{} & \multicolumn{1}{l}{} \\
\multicolumn{1}{l}{} & \multicolumn{1}{l}{} & \multicolumn{1}{l}{} &
\multicolumn{1}{l}{} & \multicolumn{1}{l}{$24\qquad ~~~n=7$} &
\multicolumn{1}{l}{} & \multicolumn{1}{l}{} \\
\multicolumn{1}{l}{$2^{(3}=4^{(3}$} & \multicolumn{1}{l}{} &
\multicolumn{1}{l}{irreducible} & \multicolumn{1}{l}{} & \multicolumn{1}{l}{$%
4\qquad \quad ~\,n=3$} & \multicolumn{1}{l}{} & \multicolumn{1}{l}{$n\quad
~~\,n\geq 3$} \\
\multicolumn{1}{l}{} & \multicolumn{1}{l}{} & \multicolumn{1}{l}{} &
\multicolumn{1}{l}{} & \multicolumn{1}{l}{$8\qquad \quad ~\,n=4$} &
\multicolumn{1}{l}{} & \multicolumn{1}{l}{} \\
\multicolumn{1}{l}{} & \multicolumn{1}{l}{} & \multicolumn{1}{l}{} &
\multicolumn{1}{l}{} & \multicolumn{1}{l}{$14\qquad ~~~n=5$} &
\multicolumn{1}{l}{} & \multicolumn{1}{l}{} \\
\multicolumn{1}{l}{} & \multicolumn{1}{l}{} & \multicolumn{1}{l}{} &
\multicolumn{1}{l}{} & \multicolumn{1}{l}{$20\qquad ~~~n=6$} &
\multicolumn{1}{l}{} & \multicolumn{1}{l}{} \\
\multicolumn{1}{l}{} & \multicolumn{1}{l}{} & \multicolumn{1}{l}{} &
\multicolumn{1}{l}{} & \multicolumn{1}{l}{$24\qquad ~~~n=7$} &
\multicolumn{1}{l}{} & \multicolumn{1}{l}{} \\
\multicolumn{1}{l}{$2^{(3}=5^{(3}$} & \multicolumn{1}{l}{} &
\multicolumn{1}{l}{reducible} & \multicolumn{1}{l}{} & \multicolumn{1}{l}{$%
2^{n-1}\qquad n\geq 3$} & \multicolumn{1}{l}{} & \multicolumn{1}{l}{$%
F_{n}\quad n\geq 3$} \\
\multicolumn{1}{l}{$3^{(3}=4^{(3}$} & \multicolumn{1}{l}{} &
\multicolumn{1}{l}{reducible} & \multicolumn{1}{l}{} & \multicolumn{1}{l}{$%
2^{n-1}\qquad n\geq 3$} & \multicolumn{1}{l}{} & \multicolumn{1}{l}{$n\quad
~~\,n\geq 3$} \\
\multicolumn{1}{l}{$3^{(3}=5^{(3}$} & \multicolumn{1}{l}{} &
\multicolumn{1}{l}{irreducible} & \multicolumn{1}{l}{} & \multicolumn{1}{l}{$%
2^{n-1}\qquad n\geq 3$} & \multicolumn{1}{l}{} & \multicolumn{1}{l}{$n\quad
~~\,n\geq 3$} \\
\multicolumn{1}{l}{$4^{(3}=5^{(3}$} & \multicolumn{1}{l}{} &
\multicolumn{1}{l}{reducible} & \multicolumn{1}{l}{} & \multicolumn{1}{l}{$%
2^{n-1}\qquad n\geq 3$} & \multicolumn{1}{l}{} & \multicolumn{1}{l}{$%
F_{n}\quad n\geq 3$}
\end{tabular}
\end{center}

\par
The behavior of $1^{(3}=5^{(3}$ shows something which wouldn't
expect - repeated checks have not revealed any error - namely
that $h_{7}=16$ is less than $h_{6}=20$. Is this an indication
that for all higher $n$ the series $h_{8},h_{9},\ldots$ is bounded
upwards? A phenomenon like that should not be excluded. For
example for the algebra defined by the two identities
$1^{(3}=4^{(3}$ and $3^{(3}=5^{(3}$ the classnumbers are
$h_{3}=3$, $h_{4}=3$ , $h_{n}=1$ for $n\geq 5$. In final analysis
the method of tableaux, as employed here, boils down to {\it set
theory.} Namely substituting two sets $M_{1},M_{2}$ of label
numbers by the single set $M_{1}\cup M_{2}-M_{1}\cap M_{2}$ {\it
provided } $M_{1}\cap M_{2}$ {\it is not empty}. Surprises are to
be expected!
\par
\smallskip
For algebras of order 4 we have conducted some investigations to
determine their classnumbers up to order 7. Since $S_{4}=14$
there are $\frac{14\times 13}{2}=91$ such identities, reducible
and irreducible. We checked 10 of them, taken at random, and
found that for 8 of them (f.ex. $11^{(4}=14^{(4}$) the
classnumbers seem to follow the rule

$$
h_{n}=S_{n}-2\frac{4^{n-3}-1}{3}+\frac{(n-2)(n-3)}{2},\qquad
3\leq n\le 7.
$$

As this expression becomes negative after $n=14$ it certainly
cannot give a general formula for $h_{n}$, unless strange
phenomena appear for $n>14$, similar to those observed for the
algebra $1^{(3}=5^{(3}$. A computer search is imperative to
decide the ambivalence.
\par
\bigskip
{\bf 4. Notes and Remarks}
\par
{\bf 4.1. Algebras with several operations of any arity.} The
process of labeling iterates and the construction of the relevant
tableaux can be used for any algebra with any finite number of
operations of any arity. The order $n$ of an iterate is the
number of operation symbols figuring in each expression. The
iterate $(a_{i}\ge 0,\quad p_{i}\ge 0)$:

$$
J\left\{
\begin{tabular}{cccc}
$a_1$-$V_1$ & $a_2$-$V_2$ & $\ldots$ & $a_k$-$V_k$ \\
$p_{1}$ & $p_{2}$ & $\ldots$ & $p_{k}$
\end{tabular}
\right\}\smallskip $$

containing the $a_{1}$-ry operation $V_{1}$, $\ p_{1}$ times, the
$a_{2}$-ry operation $V_{2}$, $p_{2}$ times, $\ldots$ has order
$$
n=p_{1}+p_{2}+\cdots +p_{k}
$$
and the number of variable places is
$$ (a_{1}-1)p_{1}+(a_{2}-1)p_{2}+...+(a_{k}-1)p_{k}+1. $$
\par
It is possible to calculate the corresponding ''Catalan'' numbers for the
nunber of iterates via generating functions as done in the classical case $%
k=1$, $a_{1}=2$. Their generating function $\phi (t)$ satisfies
the equation
$$ \frac{\phi -1}{t}=\phi ^{a_{1}}+\phi^{a_{2}}+...+\phi ^{a_{k}}
$$
which reduces for $k=1$, $a_{1}=2$ to the classical equation
$\frac{\phi -1}{t}=\phi ^{2}$, giving the Catalan numbers
$S_{n}$, and similarly to the equations $\frac{\phi -1}{t}=\phi
^{a},$ for the higher ''Catalan'' numbers:

$$ S_{n}^{a}=\frac{1}{(a-1)n+1}{an \choose n},\qquad a\ge 3.$$

\par
The proof of above functional equation for the generating function is
basically the same as in the case of the ordinary Catalan numbers, but
explicit formulas for the coefficients are only possible in exceptional
cases (see Berndt [1965], p. 71).
\par
We believe that the method of tableaux used in section 2 to
calculate the number of reducible identities can be extended to
the general case. Will the end result be the same, namely that
the number of irreducible identities divided by the number of
possible identities goes to zero with increasing $n$ ? We do not
know. But if such is the case we are confronted again with a
limitation in postulating new algebras of some order without
risking to fall back, formally or semantically, on lower order
algebras already established as essential. A preliminary check in
the case of two binary operations $V_{1}x_{1}x_{2}$ and
$V_{2}x_{1}x_{2}$ seems to confirm this contention.
\par
\smallskip
{\bf 4.2. Tableaux and Formal Languages.} The method of tableaux
we used to prove the results of sections 2 and 3 is a special case
of a general method applicable to formal languages.
\par
Roughly described, it consists in assigning a ''length'' $n$ to
the formulas ($n$-formulas) of the system and label them
recursively from $n$ to $n+1$ in a specific order proper to the
system. A formal deduction rule operating in the system generates
$(n+1)$-formulas from $n$-formulas. This allows the recursive
formation of tableaux by listing the labels of the $n$-formulas
in a line, in increasing order, and beneath it the labels of the
corresponding deducible $(n+1)$-formulas. A $(n+1)$-formula may
be deducible from several $n$-formulas and thus have its label
appear in several places of the tableau. The tableaux may be
square, orthogonal, of unequal column lengths, finite or
infinite, depending on the system. An equivalence relation
defined over the Cartesian product of the set of $n$-formulas,
subject to the rule that its validity extends to its validity for
the corresponding deducible $(n+1)$-formulas, partitions all
pairs of $n$-formulas into two classes. One class will comprise
the pairs which imply the validity of the equivalence for {\it at
least one pair} of $(n-1)$-formulas. The other class will
comprise the pairs which do not have this property. This
introduces a concept of {\it reducibility} of a formal system
equipped with an equivalence relation, to a system of lower
order. In the case of finite tableaux we can consequently count
the number of pairs in either class.
\par
Above procedure applied to the formal language of expressions $x$,
$Vxx$, $VVxxx$, $VxVxx$, $VVVxxxx$, $\ldots$ generated by
iterating the binary operation $Vxy$ yields the tableaux of our
paper and leads to the asymptotic formula for the number of
formally irreducible algebras and the other results, the
equivalence relation in question being the identity between two
$n$ -iterates of the operation.
\par
\smallskip
{\bf 4.3. Commutativity.} In section 3 we have investigated
$n$-iterates where the variable spaces are filled respectively by
the same variable on both sides of the identity. But what happens
if for example we have a structure (algebra) where following
holds:
$$
V(x_{1},VVx_{2}x_{3},x_{4})=V(x_{3},Vx_{2},Vx_{4}x_{1})
$$
identically in all variables? The simplest case is of course
commutativity $Vx_{1}x_{2}=Vx_{2}x_{1}$. The general case is
$$
J^{n}(x_{1},\ldots,x_{n+1})=J^{n}(P^{n}(x_{1}),\ldots,P^{n}(x_{n+1}))
$$
where $P^{n}$ is a permutation of the symmetric $(n+1)$-group of
the variables $x_{1},x_{2}\dots,x_{n+1}$. The generation of
$(n+1)$-iterates from $n$-iterates is the same as in 2.1 but
relabeling of the variables has to be made via following
recursion in order to get the $k$-th line of the corresponding
tableau $A_{n+1}$:
$$ P_{k}^{n+1}(x_{\nu })=\left\{
\begin{tabular}{llc}
$P^{n}(x_{\nu })$ & if & $1\le P^{n}(x_{\nu })\le k$ \\
$P^{n}(x_{\nu +1})$ & if & $k+1\le P^{n}(x_{\nu +1})\le n$ \\
$x_{k+1}$ & if & $\nu =n+1$
\end{tabular}
\right.
$$

If $P^{n}(x_{\nu })=x_{\nu }$ we get the re-labeling of the
variables used in 2.1.
\par
As no single permutation can be excluded the total number of
possible $n$ -iterates is now $(n+1)!S_{n}$ and the relevant
$A_{n}$ tableau has $n!S_{n-1}$ lines and $n$ columns.
\par
Whether an asymptotic formula holds as in 2.6 is an open question
of great interest in view of the importance of commutativity and
its higher analogues as well as other identities such as the
Jacobi identity with two operations.
\par
\smallskip
{\bf 4.4. Examples.} Suppose a structure $S$ formally defined as
follows: If $x\in S$ and $y\in S$ then $Axy\in S$ where $A$ is a
fixed symbol. Putting $Vxy=Axy $ we get the basic system of
iterates with $V$ replaced by $A$. If $L(x)$ is the number of
$A$'s figuring in $x$ we obviously have
$$
L(Vxy)=L(x)+L(y)+1
$$
Setting both sides equal to $n$ we get successively:
$$L(Vxy)=L(x)+L(y)+1=n  $$
$${\sum_{\{z|L(z)=n\}}
 }R\left[ Vxy=z\right] ={\sum_{\kappa +\lambda +1=n} }R\left[ L(x)=\kappa \right] R\left[ L(y)=\lambda \right]
$$
where $R\left[ E(x,y,...)=a\right] $ denotes generally the number
of solutions $\left[ x_{i},y_{i}...\right] $ of any equation
$E(x,y,...)=a$ with a fixed $a$. Since in the present case
$R\left[ Vxy=z\right] =1$ - an iterate $J$ can be formally
expressed in only one way as $V(J_{1},J_{2})$- we get
$$
\sum_{\left\{z|L(z)=n\right\}}1={\sum_{\kappa +\lambda +1=n} }R
\left[ L(x)=\kappa \right] R\left[ L(y)=\lambda \right] ,
$$
 The left hand side being $R\left[ L(z)=n\right] $, i.e. the $n-$th
Catalan number $S_{n}$, we get
$$
S_{n}={\sum_{\kappa +\lambda +1=n} }S_{\kappa }S_{\lambda }
$$
which is the classical recursion for these numbers.
\par
\smallskip
Above reasoning and results can be widely generalized and in our
view may be useful in investigating sequences over an alphabet
generated by a law. Let us consider for example a structure
defined by a binary law $Vxy$ and an application $L(x)$ onto the
natural numbers subject to following homomorphism
$$
L(V(xy)=P(L(x),L(y)),
$$
where $P(s,t)$ is a polynomial in $s$ and $t$ over the ring of
integers satisfying a suitable positivity condition. $P$ need not
be symmetric in $x$, $y$ or have only positive coefficients. But
e.g. let us take $P(s,t)=s^{2}+t^{2}+2$. A structure $S$ obeying
above homomorphism can be constructed as follows:
\par
\begin{center}
\begin{tabular}{ll}
1.\quad & The elements $x,y,\ldots\in S$ are sequences of $a$'s
and $A$'s. \\ & $a$ is the generator, $A$ a fixed symbol. \\
2.\quad & $L(x)$ is the number of $A$'s figuring in $x$. \\
3.\quad & If $x\in S$ and $y\in S$, $L(x)=s,L(y)=t,$ then \\ &
$V(x,y)=A\underbrace{x \cdots x}_{s times}A \underbrace{y \cdots
y}_{t times}\in S$
\end{tabular}
\end{center}
It is clear that $L(Vxy)=P(L(x),L(y))=L(x)^{2}+L(y)^{2}+2$ \ as
required.
\par
Permutation of the $A$'s in the string of symbols $x$ and $y$ give
rise to different structures satisfying
homomorphisms with the same polynomial law $P(s,t)$. In the classical case $%
Vxy=Axy$ generates parenthesizing. But $\overline{V}xy=xAy$ does not.\bigskip
\par
\smallskip
{\bf 4.5. ''Catalan'' numbers relative to homomorphisms.} The
method used to obtain the recursion formula for the Catalan
numbers can be applied to the general case, provided every
element $a$ of the structure can be uniquely expressed as
$V(x_{1}x_{2})$. In other words the equation $V(u,v)=a$ should
have only one solution $[ x_{1},x_{2}]:V(x_1x_2)=a $.
\par
Starting with a homomorphism $L$ and a polynomial $P$, and
proceeding as before we set
$$
L(Vxy)=P(L(x),L(y))=n
$$
and get the ''counting'' equation
$$
\sum_{\{z,|L(z)=n\}}R\left[ Vxy=z\right] =
\sum_{\{s,t|P(s,t)=n\}}R\left[ L(x)=s\right] R\left[ L(y)=t\right]
$$
where $R[ \ldots = \ldots ] $ is the number of solutions of the
equation in parenthesis and $\{\ldots|\ldots \}$ the all -
quantifier. Since we assumed that $R[Vxy=z]=1$ the left side
reduces to
$$
\sum_{\{z|L(z)=n\}}1=R[L(z)=n]
$$
which we denote by $S_{n}^{L}$. We suggest to call these numbers
the ''Catalan'' numbers of the structure relative to the homomorphism $(L,P)$%
. The recurrence relation now becomes\medskip
$$
 S_{n}^{L}=\sum_{P(s,t)=n}S_{s}^{L}S_{t}^{L},
$$
where the summation extends over the solutions $\geq 0$ of the
Diophantine equation $P(s,t)=n$. If $P(s,t)=s+t+1$ we get the
ordinary Catalan numbers.
\par
By way of generators $a,b,\cdots$ and fixed parameters (symbols)
$A,B,...$ we may construct algebras whose elements are words
formed with the generators and the parameters. $L(x)$ may count
the number of $A$'s in $x$, or the number of $B$'s, or for that
matter the number of any configuration of $A$'s and $B$'s, such
as for example $AB$. The composition law $Vxy$ can be set in such
a way that given a polynomial $P(s,t)$, the homomorphism $
L(Vxy)=P(Lx,Ly) $ is fulfilled, as done in 4.4. In the general
case $R[Vxy=z]$ for $L(z)=n$ is a function of both $z$ and $n$ not
always equal to $1$, depending on the underlying structure.
Recursion formulas and generating functions, in the form of power
or Dirichlet series for the ''Catalan'' numbers
$S_{n}^{L}=R[L(x)=n]$ are not always deducible, but in certain
cases can be found. If e.g. $P(s,t)=s+t+k$ and $R[Vxy=z]=n+l$,
with fixed integers $k$ and $l$, and, assuming such a structure
can be constructed as outlined, the counting equation leads to
following differential equation for the generating function
$G(t)=\sum S_{n}^{L}t^{n}$
$$
tG^{\prime}(t)+lG(t)-t^{k}G(t)^{2}=\sum\limits^{k-1}_{n=0}(n+l)S_{n}^{L}t^{n}.
$$
\par
\smallskip
{\bf 4.6. Sets of identities defining an algebra with a binary
operation.} Suppose an $n$-algebra is defined by classes $h_{n}$
of semantically equal $n$-iterates each class containing $\lambda
_{i}$ different iterates and each of the $S_{n}$ iterates
figuring in one class only. Then
$$ \lambda _{1}+\lambda _{2}+\cdots+\lambda _{h_{n}}=S_n.$$
\par
Suppose further the operation $Vx_{1}x_{2}$ satisfies the unicity
condition for the solution $Vxa=b$ and form the $S_{n-1}$
$n$-iterates $V(J_{1}^{n-1},\,x),...,V(J_{s_{n-1}}^{n-1},\,x)$.
These iterates are semantically different as otherwise there
would be at least a pair $J_{i}^{n-1}=J_{j}^{n-1}$ and the
algebra would be of lower order. Hence
$$
h_{n}\ge S_{n-1}=\frac{n+2}{4n+2}S_{n}. $$

Substituting for $S_{n}$ its above value we successively get
$$
h_{n}\min (\lambda _{i})\,\leq \,\lambda _{1}+\lambda
_{2}+...+\lambda _{h_{n}}\,\leq \,\frac{4n+2}{n+2}h_n$$

$$ \min (\lambda _{i})\leq \frac{4n+2}{n+2}<4.  $$

 Therefore $\min (\lambda _{i})\le 3, $ i.e. the set of classes defining the
algebra must contain at least one class containing max. 3
iterates, if unicity conditions apply for the binary operation.
\par
\smallskip
{\bf 4.7. Multi-variable homomorphism.} Homomorphisms $L(x)$ are
1-variable functions. We were unable to find any reference in the
literature to a generalization to more variables. In our view a
correct extension $L^n(x_1,\ldots,x_n)$ to several variables of
this fundamental concept would be as follows:
$$
L^{n}(V^{m}(x_{1}y_{1}\ldots w_{1}),\ldots,(V^{m}(x_{n}y_{n}\ldots
w_{n}))=W^{m}(L^{n}(x_{1}x_{2}\ldots x_{n}),\ldots
,L^{n}(w_{1}w_{2}\ldots w_{n})),
$$
where $V^m$ and $W^m$ are $m$-ry operations, the upper indices
indicating the number of variables of the operations involved.
This definition is purely formal, without regard to the domains of
validity of the operations, which have to be prescribed in a
concrete application in order to make it an identity.
\par
For $n=1$, $m=2$ we get the usual homo $(h)$-morphism
$L(Vxy)=W(Lx,Ly)$ for binary structures. There is however a
significant distinction between $n=1$ and $n\geq 2$, if $V=W$ so
that all functions operate in the same domain. If $n=1$ the
series $L(x),L(L(x)),...L^{n}(x)=L(L^{n-1}(x)),...$ i.e. the group
of auto $(a)$-morphisms never leads outside the set of 1-variable $a$%
-morphisms. However, this is not true for $n\geq 2$ where from the
existence of one $n$-variable $a$-morphism we can deduce new
$a$-morphisms
of any order. If e.g. $L(x_{1}x_{2})$ is a 2-variable $a$-morphism then $%
\overline{L}(x_{2}x_{2}x_{3})=L(L(x_{1}x_{2}),x_{3})$ is a
3-variable $a$-morphism. The general result establishing the
existence of a tower of $a$-morphisms is a consequence of the
definition. In the case $V=W$ the
definition is reciprocal in the sense that if $L$ is a $n$-variable $a$%
-morphism then $V$ is a $m$-variable $a$-morphism of $L$.
\par
The search of $h$-morphisms resp. $a$-morphisms with many variables of a
structure and the investigation of the corresponding groups by way of their
composition (iteration) is of importance as shown by the classical case in
the Galois Theory of equations over a field. There are analogous theorems if
$n\geq 2$, $m\geq 2$:
\par
1) Repeated application of the definition with $V=W$ \ leads to the fact
that $L$ is also an $a$-morphism of any iterate $J$ of $V$ which in short we
write $L(J)=J(L)$.
\par
2) In an equation between iterates of $V$, which may also involve constants $%
a,b,c,l\dots$ as coefficients, if $r_1,r_2,\ldots$ are roots of the equation then $%
L(r_1,r_2,\ldots)$ are roots of the same equation whose
corresponding coefficients are
$L(a,a,\ldots),L(b,b,\ldots)\ldots$ .
\par
3) If $L(a,a,\ldots)=a$ (idempotent) then $L(r_1,r_2,\ldots)$ is
also a root of the same equation.
\par
4) Considering the roots of an equation as functions of its coefficients it
can be shown that $L$ is an $a$-morphism of these functions and vice-versa,
if unicity conditions are fulfilled.
\par
Above statements are only indicatively outlined but may convey an
idea of what can be done with multi-variable homomorphisms which
of course are identities with several operations as described in
4.1. We think they can prove useful for finite structures. We
glimpse an eventual connection with $n$-Categories (see Baez,
Dolan, [1998]).
\par
\smallskip
{\bf 4.8. Relations between Catalan numbers.} We give hereunder
without proof some relations for the classical Catalan numbers
$S_{n}=\frac{1}{n+1}{2n \choose n}$, $n=1,2,...$ which we haven't
seen in the literature we are acquainted with:
$$
1.\qquad\qquad \sum_{i_{1}+\cdots+i_{\lambda+1 }=n-\lambda }
S_{i_{1}}S_{i_{2}}\ldots S_{i_{\lambda +1}}= S_{n}-{\lambda+1
\choose 1} S_{n-1}+\cdots+{(-1)^k}{{\lambda -k} \choose
k}S_{n-k}\qquad\qquad
$$
\qquad  where $k=\left[ \frac{\lambda }{2}\right] $, the greatest
integer $\le \frac{\lambda }{2}$ and $1\leq \lambda \leq n$.
\par
\qquad For $\lambda =1$ this reduces to  $\sum_{i_{1}+.i_{2}=n-1}
S_{i_{1}}S_{i_{2}}=S_{n}.$

$$ 2.\qquad \qquad \qquad S_{n}=c_{\lambda ,1} \sum_{i_{1}+i_{2}=n-\lambda }
S_{i_{1}}S_{i_{2}}+\cdots +c_{\lambda ,\lambda
}\sum_{i_{1}+\cdots +i_{\lambda +1}=n-\lambda }
S_{i_{1}}S_{i_{2}}\ldots S_{i_{\lambda +1}}\qquad \qquad \qquad $$

\qquad where $c_{\lambda ,i}$ are the numbers of the ''Pascal''
triangle defined in 2.2. \par \qquad
For $\lambda =1$ this again
reduces to the definition of the Catalan numbers. \medskip

$$
3.\qquad\qquad\qquad\qquad\lim_{n\to\infty }\frac{1}{S_{n}}
\sum_{i_{1}+\cdots +i_{\lambda +1}=n-\lambda
}S_{i_{1}}S_{i_{2}}\ldots S_{i_{ \lambda +1}}=\frac{\lambda
+1}{2^{\lambda }}\quad,\qquad \quad\lambda \ge
1.\qquad\qquad\qquad\qquad
$$

4.\quad An iterate of $Vxx$ of order $n$ is a string of $2n+1$
symbols each of which \par \qquad is either a $V$ or an $x.$ It
can therefore be written as
\begin{center}
 $ a_{1}b_{1}\cdots a_{k}b_{k}$\qquad with integers
$a_{i},b_{i} \ge 1$,
\end{center}
\par\qquad thus uniquely assigning an integer $1 \le k < n $ to each
iterate.
\par\qquad
For example $J_{11}^{(5}=VVVxxVxVxxx$ can be written as
$a_{1}b_{1}a_{2}b_{2}a_{3}b_{3}=321113$
\par\qquad with $k = 3$.

\par \qquad An expression $a_{1}b_{1}\ldots a_{k}b_{k}$ represents an
$n$-iterate of a binary operation,
\par \qquad is well formed so to say, iff
$k$ and the $a_{i},b_{i}$ satisfy one of following conditions:
\par
\begin{center}
 i) $k=1, \qquad a_1=n,\quad b_1=n+1$ \qquad\qquad\qquad
\end{center}
\newpage
\begin{center}
  ii) $k \ge 2$,\qquad\qquad\qquad $a_1+a_2+\cdots +a_k=n$ \qquad\qquad\qquad\qquad\qquad\qquad

$b_1+b_2+\cdots +b_k=n+1$

$a_1\ge b_1$

$a_1+a_2\ge b_1+b_2$

..............................

$a_{1}+a_{2}+\cdots +a_{k-1}\ge b_{1}+b_{2}+\cdots +b_{k-1}$
\end{center}

\par \qquad Hence $S_{n}-1$ is the number of integer solutions $\geq 1$ of above
Diophantine system,\par \qquad which simultaneously gives an
algorithm to decide in $k$ steps whether any word \par \qquad in
$V$'s and $x$'s is an $n$-iterate of the binary operation $Vxx$
the number of which \par \qquad is $S_n$ whose order is $\sim
{\frac {1}{\sqrt
 \pi}}n^{-3/2}2^{2n}$.
\par
\smallskip
{\bf 4.9. Skein polynomials of iterates of a binary operation.}
The build-up of an iterate $J$ can be characterized by a
polynomial $P(J;a,b)$ of 2 variables in the following way. E.g.
let us take again $J_{11}^{(5}$ and put indexes in the $V$'s and
$x$'s

$$
V_{1}V_{2}V_{3}x_{1}x_{2}\;V_{4}x_{3}V_{5}x_{4}x_{5}\;x_{6}.
$$

The relation of a variable to the operations is not always the same. For
example, $x_{2}$ is within the ''Wirkungsbereich'' (WB) of $V_{1},V_{2}$ and
$V_{3}$ but outside of that of $V_{4}$ and $V_{5}$. Similarly $x_{6}$ is
within $V_{1}$ but outside of $V_{2},V_{3},V_{4},V_{5}$. This can easily
been seen if the parentheses are put in place. A further distinction to be
made is whether the variable is in the first place or the second place of an
operation. In above example $x_{3}$ is in the first place of $V_{1}$, in the
second place of $V_{2}$ and in the first place of $V_{4}$.
\par
We substantiate these remarks by associating to each variable a
monomial $a^{s}b^{t}$ with $a,b$ in a commutative ring $R$. We
could also assign monomials $a^{s}b^{t}c^{k}$, $c$ standing for
being outside of the WB but this doesn't make any essential
difference. The exponents $s$ and $t$ are the numbers of the
first resp. second place occurrences of the variable. We then sum
these monomials to form a polynomial corresponding to the iterate
in question. Setting 1 for the polynomial of $x$ we get by this
procedure following table for the iterates $J_{i}^{n}$ up to
$n=3$, listed in the ordering defined in 2.1 (Exhibit 2):
\newpage
\begin{center}
\begin{tabular}{lll}
Iterate $J$ &  & Polynomial $P(J)$ \\
------------ &  & --------------------- \\
$x$ &  & $1$ \\
&  &  \\
$Vxx$ &  & $a+b$ \\
&  &  \\
$VVxxx$ &  & $a^{2}+ab+b$ \\
$VxVxx$ &  & $a+ab+b^{2}$ \\
&  &  \\
$VVVxxxx$ &  & $a^{3}+a^{2}b+ab+b$ \\
$VVxxVxx$ &  & $a^{2}+2ab+b^{2}$ \\
$VVxVxxx$ &  & $a^{2}+a^{2}b+ab^{2}+b$ \\
$VxVVxxx$ &  & $a+a^{2}b+ab^{2}+b^{2}$ \\
$VxVxVxx$ &  & $a+ab+ab^{2}+b^{3}$ \\
$\cdots$ &  & $\cdots$
\end{tabular}
\end{center}

It is easily seen that the formation rule of $P(J)$ obeys the
homomorphism
$$ P(V(J,J^\prime ))=aP(J)+bP(J^\prime ).  $$
\par
If therefore $R$ is the ring of integers the arguments of 4.3 can
be applied verbatim. Writing $P(J_{i}^{n};a,b)$ for the
polynomial of $J_{i}^{n}$ we have $P(J_{i}^{n};1,1)=n+1$ and
$P(J_{i}^{n};a,b)=1$ if $a+b=1$. Because of this, simpler
polynomials $Q(J)$ can be defined by setting
$$
Q(J_{i}^{n};a,b)=\frac{P(J_{i}^{n};a,b)-1}{a+b-1}
$$
with
$$
Q(V(J,J^{^{\prime }}))=aQ(J)+bQ(J^{^{\prime }})+1
$$
The function
$$
G_{P}(t)=\sum_{J}t^{P(J)}=t^{1}+t^{a+b}+t^{a^{2}+ab+b}+t^{a+ab+b^{2}}+\cdots
$$
satisfies the functional equation
$$ G_{P}(t^{a})G_{P}(t^{b})=G_{P}(t)-t.  $$
The function
$$
G_{Q}(t)=\sum_{J}t^{Q(J)}=1+t+t^{a+1}+t^{b+1}+\cdots
$$
satisfies the functional equation
$$ tG_{Q}(t^{a})G_{Q}(t^{b})=G_{Q}(t)-1,  $$
which reduces to the Catalan equation for $a=b=1$.
\par
Collecting terms in $G_{P}(t)$ for which $P(J)=P(J^{^{\prime }})$
we may write
$$ G_{P}(t)=\sum_{J}N_{P}\,t^{P(J)}   $$
where $N_{P}$ is the number of solutions $J$ for which
$P(J)=P_{i}$ , where $P_{i}$ is a specific polynomial of the
series $P_{1},P_{2},\cdots$ listed above. The $N_{P}$ satisfy the
equation:
$$ N_{P}=\sum_{\left\{ \kappa ,\lambda \,|\;aP_{\kappa
}+bP_\lambda =P\right\} }N_{P_{\kappa }}N_{P_\lambda }.
$$
\par
It was originally thought that there was a (1, 1) correspondence
between the $J$'s and their polynomials $P(J)$ and our efforts to
find the essential algebras $J_{i}^{n}=J_{j}^{n}$ consisted in
equating the corresponding polynomials
$P(J_{i}^{n};a,b)=P(J_{j}^{n};a,b)$ and look for distinctive
properties of the resulting algebraic manifolds. The
correspondence however is not (1, 1) which would mean that all
$N_{i}=1$. Already for $n=4$,
$$
P(J_{4}^{4})=P(J_{7}^{4})=a^{2}+ab+a^{2}b+ab^{2}+b^{2}\quad
(=P^\prime)
$$
so that $N_{P^\prime}=2$.
\par
Our attention has been drawn to the fact that the polynomials we
constructed to reflect the way the iterates of a binary operation
are built, bear a similarity with the polynomials of knot theory.
\par
\smallskip
{\bf 4.10. Which algebras are essential? A conjecture.} The
theorem of 3.2 says that the only essential algebras of order 3,
subject to unicity conditions, are the following two:
$$
V(V(Vx_{1}x_{2},x_{3}),\,x_{4})=V(x_{1},\,V(Vx_{2}x_{3},x_{4}))\qquad
(1^{(3}=4^{(3})
$$
$$
V(V(x_{1},\,Vx_{2}x_{3}),x_{4})=V(x_{1},\,V(x_{2},\,Vx_{3}x_{4}))\qquad
(3^{(3}=5^{(3}).
$$
We notice that as to their form these can be written respectively
as
$$
V(A,x)=V(x,A)
$$
$$
 V(B,x)=V(x,B),
$$
where $A$ and $B$ are the two iterates of order 2,  $J_{1}^{2}$
and $J_{2}^{2}$. On the other hand it is seen that the pairs (1,
4) and (3, 5) are the columns of tableau $B_{3}$. Is this a lead
sign that even in the general case the essential algebras are
given by the columns of tableaux $B_{n}$ ? \ So that for $n=4$
the essential algebras would be $J_1^4=J_9^4,\, $
$J_3^4=J_{10}^4,\, $ $J_6^4=J_{12}^4,\, $ $J_8^4=J_{13}^4,\, $
$J_{11}^4=J_{14}^4,\, $ and for any $n$ they are given by the
$S_{n-1}$ pairs\medskip
$$
V(J_{i}^{n-1},\,x\,)=V(x,\,J_{i}^{n-1}\,)\qquad i=1,...S_{n-1}.
$$
We suspect this to be true. Use of computer should help at least
for low $n$'s.
\par
\bigskip
\bigskip
\newpage
{\bf References.}
\par
\smallskip
\begin{tabular}{ll}
1. & Aigner, M. [2001]: Catalan and other numbers: a recurrent
theme, in Algebraic \\ & Combinatorics and Computer Science,  Eds.
H. Crapo, D. Senato, \\ & Springer Verlag Italia, Milano
\\ & \\
2. & Baez, J., Dolan, J., [1998]: Higher - dimensional algebra III: $n$%
-Categories and the \\ & algebra of opetopes, Adv. Math. 135
(1998), 145-206 \\ & \\
3. & Bakhturin Yu. R., Ol'shanskij A. Yu. [1991]:Identities,
Encyclopaedia of Mathematical \\ & Sciences Vol 18, Algebra II,
Springer - Verlag, Berlin Heidleberg
\\ & \\
4. & Berndt, B. C. [1985]: Ramanujan's Notebooks, Part I, Springer
- Verlag, New York
\\ & \\
5. & Cohn, P. M., [1965]: Umiversal Algebra, Harper\&Row, Ltd.,
London
\\ & \\
6. & Riordan, J. [1968]: Combinatorial Identities, J, Wiley, New
York
\\ & \\
7. & Tarski, A. [1968]: Equational Logic and Equational Theories
of Algebras, in \\ & Contributions to Mathematical Logic, Studies
in Logic and the Foundations of \\ & Mathematics, North - Holland
Publishing Co., Amsterdam
\\ & \\
8. & Wilf, H. [1990]. Generatingfunctionology, Academic Press,
Inc., San Diego

\end{tabular}

\newpage

\begin{center}
 Exhibit 1
\par
  {\bf Ordering and labels of 4-iterates generated
from ordering of 3-iterates.}
\par
\bigskip
\bigskip
\bigskip

\begin{tabular}{ccccc}
1 $VVVxxxx$ & 2 $VVxxVxx$ & 3 $VVxVxxx$ & 4 $VxVVxxx$ & 5 $VxVxVxx$ \\
-------------------- & -------------------- &
-------------------- & -------------------- &
-------------------- \\ &  &  &  &   \\
\multicolumn{1}{l}{1 $VVVVxxxxx$} & \multicolumn{1}{l}{~~2 $VVVxxxVxx$} &
\multicolumn{1}{l}{~~3 $VVVxxVxxx$} & \multicolumn{1}{l}{~~4 $VVxxVVxxx$} &
\multicolumn{1}{l}{~~5 $VVxxVxVxx$} \\
\multicolumn{1}{l}{6 $VVVxVxxxx$} & \multicolumn{1}{l}{~~7 $VVxVxxVxx$} &
\multicolumn{1}{l}{~~8 $VVxVVxxxx$} & \multicolumn{1}{l}{~~9 $VxVVVxxxx$} &
\multicolumn{1}{l}{10 $VxVVxxVxx$} \\
\multicolumn{1}{l}{3 $VVVxxVxxx$} & \multicolumn{1}{l}{~~4 $VVxxVVxxx$} &
\multicolumn{1}{l}{11 $VVxVxVxxx$} & \multicolumn{1}{l}{12 $VxVVxVxxx$} &
\multicolumn{1}{l}{13 $VxVxVVxxx$} \\
\multicolumn{1}{l}{2 $VVVxxxVxx$} & \multicolumn{1}{l}{~~5 $VVxxVxVxx$} &
\multicolumn{1}{l}{~~7 $VVxVxxVxx$} & \multicolumn{1}{l}{10 $VxVVxxVxx$} &
\multicolumn{1}{l}{14 $VxVxVxVxx$} \\
\multicolumn{1}{l}{} & \multicolumn{1}{l}{} & \multicolumn{1}{l}{} &
\multicolumn{1}{l}{} & \multicolumn{1}{l}{} \\
\multicolumn{1}{l}{1 $VVVVxxxxx$} & \multicolumn{1}{l}{~~3 $VVVxxVxxx$} &
\multicolumn{1}{l}{~~6 $VVVxVxxxx$} & \multicolumn{1}{l}{~~8 $VVxVVxxxx$} &
\multicolumn{1}{l}{11 $VVxVxVxxx$} \\
\multicolumn{1}{l}{9 $VxVVVxxxx$} & \multicolumn{1}{l}{10 $VxVVxxVxx$} &
\multicolumn{1}{l}{12 $VxVVxVxxxx$} & \multicolumn{1}{l}{13 $VxVxVVxxx$} &
\multicolumn{1}{l}{14 $VxVxVxVxx$}
\end{tabular}

\par
\bigskip
\bigskip
\bigskip

\begin{tabular}{rrrrr}
\multicolumn{5}{c}{Tableau $A_{4}$} \\
\multicolumn{5}{c}{---------------} \\
\textbf{1} & \textbf{2} & \textbf{3} & \textbf{4} & \textbf{5} \\
\textbf{6} & \textbf{7} & \textbf{8} & \textbf{9} & \textbf{10} \\
3 & 4 & \textbf{11} & \textbf{12} & \textbf{13} \\
2 & 5 & 7 & 10 & \textbf{14}
\end{tabular}

\par
\bigskip
\bigskip

\begin{tabular}{rrrrr}
\multicolumn{5}{c}{Tableau $B_{4}$} \\
\multicolumn{5}{c}{---------------} \\
1 & 3 & 6 & 8 & 11 \\
9 & 10 & 12 & 13 & 14 \\
&  &  &  &  \\
&  &  &  &
\end{tabular}
\end{center}

\newpage

\begin{center}
 Exhibit 2
\end{center}
\par
\bigskip

\begin{center}
\begin{tabular}{llllll}
\multicolumn{2}{l}{\textbf{Iterates} $J_{i}^{n}\quad (n\leq 4)$} &  &  &  &
\textbf{Labels} $i^{(n}$ \\
&  &  &  &  &  \\
$J_{1}^{0}$ & $x$ &  &  &  & \multicolumn{1}{c}{$~~1^{(0}$} \\
&  &  &  &  & \multicolumn{1}{c}{} \\
$J_{1}^{1}$ & $Vxx$ &  &  &  & \multicolumn{1}{c}{$~~1^{(1}$} \\
&  &  &  &  & \multicolumn{1}{c}{} \\
$J_{1}^{2}$ & $VVxxx$ &  &  &  & \multicolumn{1}{c}{$~~1^{(2}$} \\
$J_{2}^{2}$ & $VxVxx$ &  &  &  & \multicolumn{1}{c}{$~~2^{(2}$} \\
&  &  &  &  & \multicolumn{1}{c}{} \\
$J_{1}^{3}$ & $VVVxxxx$ &  &  &  & \multicolumn{1}{c}{$~~1^{(3}$} \\
$J_{2}^{3}$ & $VVxxVxx$ &  &  &  & \multicolumn{1}{c}{$~~2^{(3}$} \\
$J_{3}^{3}$ & $VVxVxxx$ &  &  &  & \multicolumn{1}{c}{$~~3^{(3}$} \\
$J_{4}^{3}$ & $VxVVxxx$ &  &  &  & \multicolumn{1}{c}{$~~4^{(3}$} \\
$J_{5}^{3}$ & $VxVxVxx$ &  &  &  & \multicolumn{1}{c}{$~~5^{(3}$} \\
&  &  &  &  & \multicolumn{1}{c}{} \\
$J_{1}^{4}$ & $VVVVxxxxx$ &  &  &  & \multicolumn{1}{c}{$~~1^{(4}$} \\
$J_{2}^{4}$ & $VVVxxxVxx$ &  &  &  & \multicolumn{1}{c}{$~~2^{(4}$} \\
$J_{3}^{4}$ & $VVVxxVxxx$ &  &  &  & \multicolumn{1}{c}{$~~3^{(4}$} \\
$J_{4}^{4}$ & $VVxxVVxxx$ &  &  &  & \multicolumn{1}{c}{$~~4^{(4}$} \\
$J_{5}^{4}$ & $VVxxVxVxx$ &  &  &  & \multicolumn{1}{c}{$~~5^{(4}$} \\
$J_{6}^{4}$ & $VVVxVxxxx$ &  &  &  & \multicolumn{1}{c}{$~~6^{(4}$} \\
$J_{7}^{4}$ & $VVxVxxVxx$ &  &  &  & \multicolumn{1}{c}{$~~7^{(4}$} \\
$J_{8}^{4}$ & $VVxVVxxxx$ &  &  &  & \multicolumn{1}{c}{$~~8^{(4}$} \\
$J_{9}^{4}$ & $VxVVVxxxx$ &  &  &  & \multicolumn{1}{c}{$~~9^{(4}$} \\
$J_{10}^{4}$ & $VxVVxxVxx$ &  &  &  & \multicolumn{1}{c}{$10^{(4}$} \\
$J_{11}^{4}$ & $VVxVxVxxx$ &  &  &  & \multicolumn{1}{c}{$11^{(4}$} \\
$J_{12}^{4}$ & $VxVVxVxxx$ &  &  &  & \multicolumn{1}{c}{$12^{(4}$} \\
$J_{13}^{4}$ & $VxVxVVxxx$ &  &  &  & \multicolumn{1}{c}{$13^{(4}$} \\
$J_{14}^{4}$ & $VxVxVxVxx$ &  &  &  & \multicolumn{1}{c}{$14^{(4}$} \\
&  &  &  &  & \\
$J_{1}^{5}$ & $VVVVVxxxxxx$ &  &  &  & \multicolumn{1}{c}{$1^{(5}$} \\
$J_{2}^{5}$ & $VVVVxxxxVxx$ &  &  &  & \multicolumn{1}{c}{$2^{(5}$} \\
$\ldots$ & \quad\qquad$\ldots$ &  &  &  & \qquad$\ldots$ \\
$J_{8}^{5}$ & $VVVxxVVxxxx$ &  &  &  & \multicolumn{1}{c}{$8^{(5}$} \\
$\ldots$ & \quad\qquad$\ldots$ &  &  &  &\qquad $\ldots$ \\
$J_{11}^{5}$ & $VVVxxVxVxxx$ &  &  &  & \multicolumn{1}{c}{$11^{(5}$} \\
$\ldots$ &\quad\qquad$\ldots$ &  &  &  & \qquad$\ldots$ \\
$J_{42}^{5}$ & $VxVxVxVxVxx$ &  &  &  & \multicolumn{1}{c}{$42^{(5}$} \\
 &  &  &  &  & \\
$\ldots$ & \quad\qquad$\ldots$ &  &  &  & \qquad$\ldots$
\end{tabular}
\end{center}

\newpage

\begin{center}
 Exhibit 3
\end{center}

\begin{center}
\begin{tabular}{ccccccccccccccc}
\multicolumn{15}{c}{\bf{Tableaux $A_{n}$ and $B_{n}\quad (n\leq 5)$}} \\
\multicolumn{15}{c}{--------------------------------------------} \\
$A_{1}$ &  &  & $A_{2}$ &  &  & \multicolumn{2}{c}{$A_{3}$} &  &  &
\multicolumn{5}{c}{$A_{4}$} \\
\multicolumn{1}{c}{\textbf{1}} & \multicolumn{1}{c}{} & \multicolumn{1}{c}{}
& \multicolumn{1}{c}{\textbf{1}} & \multicolumn{1}{c}{} & \multicolumn{1}{c}{
} & \textbf{1} & \textbf{2} & \multicolumn{1}{c}{} & \multicolumn{1}{c}{} &
\textbf{1} & \textbf{2} & \textbf{3} & \textbf{4} & \textbf{5} \\
\multicolumn{1}{c}{} & \multicolumn{1}{c}{} & \multicolumn{1}{c}{} &
\multicolumn{1}{c}{\textbf{2}} & \multicolumn{1}{c}{} & \multicolumn{1}{c}{}
& \textbf{3} & \textbf{4} & \multicolumn{1}{c}{} & \multicolumn{1}{c}{} &
\textbf{6} & \textbf{7} & \textbf{8} & \textbf{9} & \textbf{10} \\
\multicolumn{1}{c}{$B_{1}$} & \multicolumn{1}{c}{} & \multicolumn{1}{c}{} &
\multicolumn{1}{c}{} & \multicolumn{1}{c}{} & \multicolumn{1}{c}{} & 2 &
\textbf{5} & \multicolumn{1}{c}{} & \multicolumn{1}{c}{} & 3 & 4 & \textbf{11%
} & \textbf{12} & \textbf{13} \\
1 &  &  & \multicolumn{1}{c}{$B_{2}$} &  &  & &  &  &  & 2 & 5 &
7 & 10 &
\textbf{14} \\
&  &  & 1 &  &  & \multicolumn{2}{c}{$B_{3}$} &  &  &  &  &  &  &  \\

&  &  & 2 &  &  & 1 & 3 &  &    &  &  & {$B_{4}$} &  & \\
&  &  &   &  &  & 4 & 5 &  &  & 1 & 3 & 6 & 8 & 11 \\
&  &  &   &  &  &   &   &  &  & 9 & 10 & 12 & 13 & 14

\end{tabular}
\end{center}
\bigskip

\begin{center}
\begin{tabular}{cccccccccccccc}
\multicolumn{14}{c}{$A_{5}$} \\
\textbf{1} & \textbf{2} & \textbf{3} & \textbf{4} & \textbf{5} & \textbf{6}
& \textbf{7} & \textbf{8} & \textbf{9} & \textbf{10} & \textbf{11} & \textbf{%
12} & \textbf{13} & \textbf{14} \\
\textbf{15} & \textbf{16} & \textbf{17} & \textbf{18} & \textbf{19} &
\textbf{20} & \textbf{21} & \textbf{22} & \textbf{23} & \textbf{24} &
\textbf{25} & \textbf{26} & \textbf{27} & \textbf{28} \\
6 & 7 & 8 & 9 & 10 & \textbf{29} & \textbf{30} & \textbf{31} & \textbf{32} &
\textbf{33} & \textbf{34} & \textbf{35} & \textbf{36} & \textbf{37} \\
3 & 4 & 11 & 12 & 13 & 17 & 18 & 25 & 26 & 27 & \textbf{38} & \textbf{39} &
\textbf{40} & \textbf{41} \\
2 & 5 & 7 & 10 & 14 & 16 & 19 & 21 & 24 & 28 & 30 & 33 & 37 & \textbf{42} \\
&  &  &  &  &  &  &  &  &  &  &  &  &  \\
\multicolumn{14}{c}{$B_{5}$} \\
1 & 3 & 6 & 8 & 11 & 15 & 17 & 20 & 22 & 25 & 29 & 31 & 34 & 38 \\
23 & 24 & 26 & 27 & 28 & 32 & 33 & 35 & 36 & 37 & 39 & 40 & 41 & 42
\end{tabular}
\end{center}
\bigskip
\begin{center}
\begin{tabular}{ccccccccccc}
\multicolumn{11}{c}{{\bf Tableaux $A_{n}\oplus B_{n}\quad (n\leq 5)$}} \\
\multicolumn{11}{c}{---------------------------------------} \\
$A_{1}\oplus B_{1}$ &  & $A_{2}\oplus B_{2}$ &  & $A_{3}\oplus B_{3}$ &  &
\multicolumn{5}{c}{$A_{4}\oplus B_{4}$} \\
1 &  & 1 &  & 1\qquad 2 &  & 1 & 2 & 3 & 4 & 5 \\
1 &  & 2 &  & 3\qquad 4 &  & 6 & 7 & 8 & 9 & 10 \\
&  & 1 &  & 2\qquad 5 &  & 3 & 4 & 11 & 12 & 13 \\
&  & 2 &  & 1\qquad 3 &  & 2 & 5 & 7 & 10 & 14 \\
&  &  &  & 4\qquad 5 &  & 1 & 3 & 6 & 8 & 11 \\
&  &  &  &  &  & 9 & 10 & 12 & 13 & 14
\end{tabular}
\end{center}
\bigskip
\begin{center}
\begin{tabular}{cccccccccccccc}
\multicolumn{14}{c}{$A_{5}\oplus B_{5}$} \\
1 & 2 & 3 & 4 & 5 & 6 & 7 & 8 & 9 & 10 & 11 & 12 & 13 & 14 \\
15 & 16 & 17 & 18 & 19 & 20 & 21 & 22 & 23 & 24 & 25 & 26 & 27 & 28 \\
6 & 7 & 8 & 9 & 10 & 29 & 30 & 31 & 32 & 33 & 34 & 35 & 36 & 37 \\
3 & 4 & 11 & 12 & 13 & 17 & 18 & 25 & 26 & 27 & 38 & 39 & 40 & 41 \\
2 & 5 & 7 & 10 & 14 & 16 & 19 & 21 & 24 & 28 & 30 & 33 & 37 & 42 \\
1 & 3 & 6 & 8 & 11 & 15 & 17 & 20 & 22 & 25 & 29 & 31 & 34 & 38 \\
23 & 24 & 26 & 27 & 28 & 32 & 33 & 35 & 36 & 37 & 39 & 40 & 41 & 42
\end{tabular}
\end{center}
\bigskip
\newpage

\begin{center}
Exhibit 4
\end{center}
\par
\bigskip
\begin{center}
\bf {''Pascal'' triangle for the numbers} $c_{n,j}$
\end{center}

\begin{center}
\begin{tabular}{ccccccccccccc}
$n\backslash j$ &  & 1 & 2 & 3 & 4 & 5
& 6 & 7 & 8 & 9 & 10  & $\ldots$\\
&  &  & &  &  &  & & &  &  &  & \\
 &  &  &  &  &  &  &  &  &  &  &  &     \\
1 &  & 1 &  &  & &  &  &  & & &  & \\
&  &  &  &  &  &  &  &  &  &  &   &   \\
2 &  &  1 & 1 &  &  &  &  & &  &
 &  & \\
&  &  &  &  &  &  &  &  &  &  &    &  \\
3 &  & 2 & 2 & 1 &  &  &  &  &  &
 &  & \\
&  &  &  &  &  &  &  &  &  &  &   &   \\
4 &  & 5 & 5 & 3 & 1 &  &  & &
 &   &  & \\
&  &  &  &  &  &  &  &  &  &  &   &   \\
5 &  & 14 & 14 & 9 & 4 & 1 & &  &  &
 &  & \\
&  &  &  &  &  &  &  &  &  &  &    &  \\
6 &  & 42 & 42 & 28 & 14 & 5 & 1 &  &  &
 &  & \\
&  &  &  &  &  &  &  &  &  &  &      \\
7 &  & 132 & 132 & 90 & 48 & 20 & 6 & 1 &  &   &  \\
&  &  &  &  &  &  &  &  &  &  &    &  \\
8 &  & 429 & 429 & 297 & 165 & 75 & 27 & 7 & 1 &
&  & \\
&  &  &  &  &  &  &  &  &  &  &    &  \\
9 &  & 1430 & 1430  & 1001 & 572 & 275 & 110 & 35 &
8 & 1 & & \\
&  &  &  &  &  &  &  &  &  &  &    &  \\
10 &  & 4862 & 4862 & 3432 & 2002 & 1001 & 429 & 154
& 44 & 9 & 1 & \\
&  &  &  &  &  &  &  &  &  &  &   &   \\
$\cdots$ &  & $\cdots$ & $\cdots$ & $\cdots$ & $\cdots$ &
$\cdots$ & $\cdots$ & $\cdots$ & $\cdots$ & $\cdots$ & $\cdots$ &
$\cdots$
\end{tabular}
\end{center}

$$c_{n,j}=c_{n-1,j-1}+c_{n-1,j}+c_{n-1,j+1}+\cdots+c_{n-1,n-1}$$

$$\sum_{j=1}^{n}c_{n,,j}=S_{n}$$

\newpage

\begin{center}
Exhibit 5
\end{center}
\par
\begin{center}
{\bf Incidence matrix for} $n=4$ {\bf relative to tableau} $A_{4}$
{\bf for all possible} \par $S_{4}^{2}$ $(=14^{2})$ {\bf
identities} $J_{i}^{(4}=J_{j}^{(4}$ {\bf and calculation of}
$I_4$ \par ($J_{i}^{(4}$ {\bf denoted by} $i$. {\bf Blank spaces
denote} $0$'s )
\end{center}
$$
\delta _{ij}=\left\{
\begin{tabular}{lll}
$1$ & if & identity reducible \\
$0$ & if & identity irreducible
\end{tabular}
\right.
$$
\begin{center}
$M(i)=$ Multiplicity of $J_{i}$
\end{center}

\[
\begin{tabular}{ccccccccccccccccllcc}
$i\backslash j$ &  & 1 & 2 & 3 & 4 & 5 & 6 & 7 & 8 & 9 & 10 & 11 & 12 & 13 &
14 & \multicolumn{2}{c}{} & $\sum_{i}1$ & $M(i)$ \\
& \multicolumn{1}{l}{} & \multicolumn{1}{l}{} & \multicolumn{1}{l}{} &
\multicolumn{1}{l}{} & \multicolumn{1}{l}{} & \multicolumn{1}{l}{} &
\multicolumn{1}{l}{} & \multicolumn{1}{l}{} & \multicolumn{1}{l}{} &
\multicolumn{1}{l}{} & \multicolumn{1}{l}{} & \multicolumn{1}{l}{} &
\multicolumn{1}{l}{} & \multicolumn{1}{l}{} & \multicolumn{1}{l}{} &
\multicolumn{2}{l}{} &  &  \\
1 & \multicolumn{1}{l}{} & \multicolumn{1}{l}{1} & \multicolumn{1}{l}{1} &
\multicolumn{1}{l}{1} & \multicolumn{1}{l}{1} & \multicolumn{1}{l}{1} &
\multicolumn{1}{l}{} & \multicolumn{1}{l}{} & \multicolumn{1}{l}{} &
\multicolumn{1}{l}{} & \multicolumn{1}{l}{} & \multicolumn{1}{l}{} &
\multicolumn{1}{l}{} & \multicolumn{1}{l}{} & \multicolumn{1}{l}{} &
\multicolumn{2}{l}{} & 5 & 1 \\
2 & \multicolumn{1}{l}{} & \multicolumn{1}{l}{1} & \multicolumn{1}{l}{1} &
\multicolumn{1}{l}{1} & \multicolumn{1}{l}{1} & \multicolumn{1}{l}{1} &
\multicolumn{1}{l}{} & \multicolumn{1}{l}{1} & \multicolumn{1}{l}{} &
\multicolumn{1}{l}{} & \multicolumn{1}{l}{1} & \multicolumn{1}{l}{} &
\multicolumn{1}{l}{} & \multicolumn{1}{l}{} & \multicolumn{1}{l}{1} &
\multicolumn{2}{l}{} & 8 & 2 \\
3 & \multicolumn{1}{l}{} & \multicolumn{1}{l}{1} & \multicolumn{1}{l}{1} &
\multicolumn{1}{l}{1} & \multicolumn{1}{l}{1} & \multicolumn{1}{l}{1} &
\multicolumn{1}{l}{} & \multicolumn{1}{l}{} & \multicolumn{1}{l}{} &
\multicolumn{1}{l}{} & \multicolumn{1}{l}{} & \multicolumn{1}{l}{1} &
\multicolumn{1}{l}{1} & \multicolumn{1}{l}{1} & \multicolumn{1}{l}{} &
\multicolumn{2}{l}{} & 8 & 2 \\
4 & \multicolumn{1}{l}{} & \multicolumn{1}{l}{1} & \multicolumn{1}{l}{1} &
\multicolumn{1}{l}{1} & \multicolumn{1}{l}{1} & \multicolumn{1}{l}{1} &
\multicolumn{1}{l}{} & \multicolumn{1}{l}{} & \multicolumn{1}{l}{} &
\multicolumn{1}{l}{} & \multicolumn{1}{l}{} & \multicolumn{1}{l}{1} &
\multicolumn{1}{l}{1} & \multicolumn{1}{l}{1} & \multicolumn{1}{l}{} &
\multicolumn{2}{l}{} & 8 & 2 \\
5 & \multicolumn{1}{l}{} & \multicolumn{1}{l}{1} & \multicolumn{1}{l}{1} &
\multicolumn{1}{l}{1} & \multicolumn{1}{l}{1} & \multicolumn{1}{l}{1} &
\multicolumn{1}{l}{} & \multicolumn{1}{l}{1} & \multicolumn{1}{l}{} &
\multicolumn{1}{l}{} & \multicolumn{1}{l}{1} & \multicolumn{1}{l}{} &
\multicolumn{1}{l}{} & \multicolumn{1}{l}{} & \multicolumn{1}{l}{1} &
\multicolumn{2}{l}{} & 8 & 2 \\
6 & \multicolumn{1}{l}{} & \multicolumn{1}{l}{} & \multicolumn{1}{l}{} &
\multicolumn{1}{l}{} & \multicolumn{1}{l}{} & \multicolumn{1}{l}{} &
\multicolumn{1}{l}{1} & \multicolumn{1}{l}{1} & \multicolumn{1}{l}{1} &
\multicolumn{1}{l}{1} & \multicolumn{1}{l}{1} & \multicolumn{1}{l}{} &
\multicolumn{1}{l}{} & \multicolumn{1}{l}{} & \multicolumn{1}{l}{} &
\multicolumn{2}{l}{} & 5 & 1 \\
7 & \multicolumn{1}{l}{} & \multicolumn{1}{l}{} & \multicolumn{1}{l}{1} &
\multicolumn{1}{l}{} & \multicolumn{1}{l}{} & \multicolumn{1}{l}{1} &
\multicolumn{1}{l}{1} & \multicolumn{1}{l}{1} & \multicolumn{1}{l}{1} &
\multicolumn{1}{l}{1} & \multicolumn{1}{l}{1} & \multicolumn{1}{l}{} &
\multicolumn{1}{l}{} & \multicolumn{1}{l}{} & \multicolumn{1}{l}{1} &
\multicolumn{2}{l}{} & 8 & 2 \\
8 & \multicolumn{1}{l}{} & \multicolumn{1}{l}{} & \multicolumn{1}{l}{} &
\multicolumn{1}{l}{} & \multicolumn{1}{l}{} & \multicolumn{1}{l}{} &
\multicolumn{1}{l}{1} & \multicolumn{1}{l}{1} & \multicolumn{1}{l}{1} &
\multicolumn{1}{l}{1} & \multicolumn{1}{l}{1} & \multicolumn{1}{l}{} &
\multicolumn{1}{l}{} & \multicolumn{1}{l}{} & \multicolumn{1}{l}{} &
\multicolumn{2}{l}{} & 5 & 1 \\
9 & \multicolumn{1}{l}{} & \multicolumn{1}{l}{} & \multicolumn{1}{l}{} &
\multicolumn{1}{l}{} & \multicolumn{1}{l}{} & \multicolumn{1}{l}{} &
\multicolumn{1}{l}{1} & \multicolumn{1}{l}{1} & \multicolumn{1}{l}{1} &
\multicolumn{1}{l}{1} & \multicolumn{1}{l}{1} & \multicolumn{1}{l}{} &
\multicolumn{1}{l}{} & \multicolumn{1}{l}{} & \multicolumn{1}{l}{} &
\multicolumn{2}{l}{} & 5 & 1 \\
10 & \multicolumn{1}{l}{} & \multicolumn{1}{l}{} & \multicolumn{1}{l}{1} &
\multicolumn{1}{l}{} & \multicolumn{1}{l}{} & \multicolumn{1}{l}{1} &
\multicolumn{1}{l}{1} & \multicolumn{1}{l}{1} & \multicolumn{1}{l}{1} &
\multicolumn{1}{l}{1} & \multicolumn{1}{l}{1} & \multicolumn{1}{l}{} &
\multicolumn{1}{l}{} & \multicolumn{1}{l}{} & \multicolumn{1}{l}{1} &
\multicolumn{2}{l}{} & 8 & 2 \\
11 & \multicolumn{1}{l}{} & \multicolumn{1}{l}{} & \multicolumn{1}{l}{} &
\multicolumn{1}{l}{1} & \multicolumn{1}{l}{1} & \multicolumn{1}{l}{} &
\multicolumn{1}{l}{} & \multicolumn{1}{l}{} & \multicolumn{1}{l}{} &
\multicolumn{1}{l}{} & \multicolumn{1}{l}{} & \multicolumn{1}{l}{1} &
\multicolumn{1}{l}{1} & \multicolumn{1}{l}{1} & \multicolumn{1}{l}{} &
\multicolumn{2}{l}{} & 5 & 1 \\
12 & \multicolumn{1}{l}{} & \multicolumn{1}{l}{} & \multicolumn{1}{l}{} &
\multicolumn{1}{l}{1} & \multicolumn{1}{l}{1} & \multicolumn{1}{l}{} &
\multicolumn{1}{l}{} & \multicolumn{1}{l}{} & \multicolumn{1}{l}{} &
\multicolumn{1}{l}{} & \multicolumn{1}{l}{} & \multicolumn{1}{l}{1} &
\multicolumn{1}{l}{1} & \multicolumn{1}{l}{1} & \multicolumn{1}{l}{} &
\multicolumn{2}{l}{} & 5 & 1 \\
13 & \multicolumn{1}{l}{} & \multicolumn{1}{l}{} & \multicolumn{1}{l}{} &
\multicolumn{1}{l}{1} & \multicolumn{1}{l}{1} & \multicolumn{1}{l}{} &
\multicolumn{1}{l}{} & \multicolumn{1}{l}{} & \multicolumn{1}{l}{} &
\multicolumn{1}{l}{} & \multicolumn{1}{l}{} & \multicolumn{1}{l}{1} &
\multicolumn{1}{l}{1} & \multicolumn{1}{l}{1} & \multicolumn{1}{l}{} &
\multicolumn{2}{l}{} & 5 & 1 \\
14 & \multicolumn{1}{l}{} & \multicolumn{1}{l}{} & \multicolumn{1}{l}{1} &
\multicolumn{1}{l}{} & \multicolumn{1}{l}{} & \multicolumn{1}{l}{1} &
\multicolumn{1}{l}{} & \multicolumn{1}{l}{1} & \multicolumn{1}{l}{} &
\multicolumn{1}{l}{} & \multicolumn{1}{l}{1} & \multicolumn{1}{l}{} &
\multicolumn{1}{l}{} & \multicolumn{1}{l}{} & \multicolumn{1}{l}{1} &
\multicolumn{2}{l}{} & 5 & 1 \\
& \multicolumn{1}{l}{} & \multicolumn{1}{l}{} & \multicolumn{1}{l}{} &
\multicolumn{1}{l}{} & \multicolumn{1}{l}{} & \multicolumn{1}{l}{} &
\multicolumn{1}{l}{} & \multicolumn{1}{l}{} & \multicolumn{1}{l}{} &
\multicolumn{1}{l}{} & \multicolumn{1}{l}{} & \multicolumn{1}{l}{} &
\multicolumn{1}{l}{} & \multicolumn{1}{l}{} & \multicolumn{1}{l}{} &
\multicolumn{2}{l}{} & ------ &  \\
& \multicolumn{1}{l}{} & \multicolumn{1}{l}{} & \multicolumn{1}{l}{} &
\multicolumn{1}{l}{} & \multicolumn{1}{l}{} & \multicolumn{1}{l}{} &
\multicolumn{1}{l}{} & \multicolumn{1}{l}{} & \multicolumn{1}{l}{} &
\multicolumn{1}{l}{} & \multicolumn{1}{l}{} & \multicolumn{1}{l}{} &
\multicolumn{1}{l}{} & \multicolumn{1}{l}{} & \multicolumn{1}{l}{} &
\multicolumn{2}{l}{$I_{4}=$} & $88$ &  \\
\multicolumn{1}{l}{} & \multicolumn{1}{l}{} & \multicolumn{1}{l}{} &
\multicolumn{1}{l}{} & \multicolumn{1}{l}{} & \multicolumn{1}{l}{} &
\multicolumn{1}{l}{} & \multicolumn{1}{l}{} & \multicolumn{1}{l}{} &
\multicolumn{1}{l}{} & \multicolumn{1}{l}{} & \multicolumn{1}{l}{} &
\multicolumn{1}{l}{} & \multicolumn{1}{l}{} & \multicolumn{1}{l}{} &
\multicolumn{1}{l}{} &  &  & \multicolumn{1}{l}{} & \multicolumn{1}{l}{}
\end{tabular}
\]

\begin{center}
\ \ \ \   Number of $J_{i}^{n}$ with multiplicity $k$ : \ \ \ \  $%
T_{n,k}=2^{n-2k+1}{n-1 \choose 2k-2}S_{k-1}$

\ \ \  For $n=4,\quad k=1$ : \quad \ $T_{4,1}= \,\, \, \,8$

\ \ \  For $n=4,\quad k=2$ : \quad \ $T_{4,2}= \,\, \, \,6$
\par
\end{center}

$$
I_n=\sum_{k=1}^{\left[\frac{n+1}{2}\right]}(-1)^{k-1}{n-k+1
\choose k}S_{n-k}^2
$$
$$
I_4=\sum_{k=1}^2(-1)^{k-1}{5-k \choose k}S_{4-k}^2={4 \choose
1}S_3^2-{3 \choose 2}S_2^2 = 88
$$
\begin{center}
($S_2=2,\quad S_3=5$)
\end{center}
\newpage

\begin{center}
Exhibit 6
\end{center}
\par

\begin{center}
{\bf 1. Semantical equality of two iterates in same column
connected by} $\updownarrow$ {\bf implies semantical equality of
two iterates of order 2 if} $V$ {\bf satisfies unicity
conditions.}
\end{center}

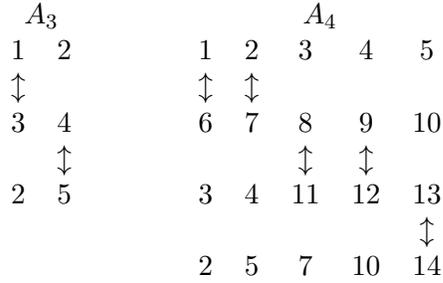
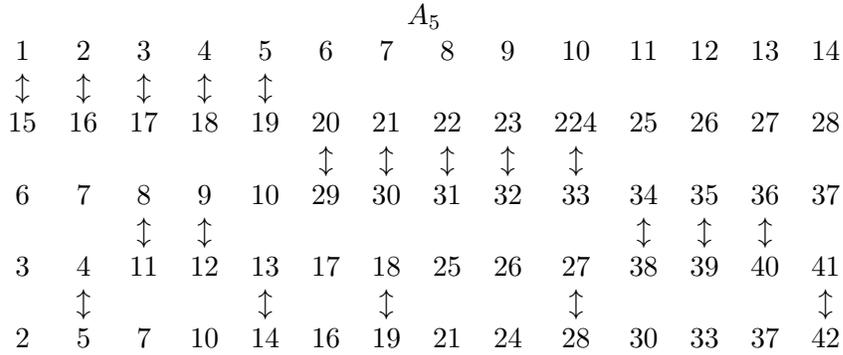
\[
\begin{tabular}{lllll}
\multicolumn{2}{c}{$A_{3}$} &  &  &    \\
1 & 2 &  &  &  \\
$\updownarrow $ &  &  &  &  \\
3 & 4 &  &  &  \\
& $\updownarrow $ &  &  &  \\
2 & 5 &  &  &  \\
&  &  &  &\\
&  &  &  &\\
&  &  &  &\\
&  &  &  &\\
\end{tabular}
\begin{tabular}{ccccc}
\multicolumn{5}{c}{$A_{4}$} \\
1 & 2 & 3 & 4 & 5 \\
$\updownarrow $ & $\updownarrow $ &  &  &  \\
6 & 7 & 8 & 9 & 10 \\
&  & $\updownarrow $ & $\updownarrow $ &  \\
3 & 4 & 11 & 12 & 13 \\
&  &  &  & $\updownarrow $ \\
2 & 5 & 7 & 10 & 14 \\
&  &  &  &  \\
&  &  &  &
\end{tabular}
\]
\par
\begin{center}
\begin{tabular}{cccccccccccccc}
\multicolumn{14}{c}{$A_{5}$} \\
1 & 2 & 3 & 4 & 5 & 6 & 7 & 8 & 9 & 10 & 11 & 12 & 13 & 14 \\
$\updownarrow $ & $\updownarrow $ & $\updownarrow $ & $\updownarrow $ & $%
\updownarrow $ &  &  &  &  &  &  &  &  &  \\
15 & 16 & 17 & 18 & 19 & 20 & 21 & 22 & 23 & 224 & 25 & 26 & 27 & 28 \\
&  &  &  &  & $\updownarrow $ & $\updownarrow $ & $\updownarrow $ & $%
\updownarrow $ & $\updownarrow $ &  &  &  &  \\
6 & 7 & 8 & 9 & 10 & 29 & 30 & 31 & 32 & 33 & 34 & 35 & 36 & 37 \\
&  & $\updownarrow $ & $\updownarrow $ &  &  &  &  &  &  &
$\updownarrow $ &
$\updownarrow $ & $\updownarrow $ &  \\
3 & 4 & 11 & 12 & 13 & 17 & 18 & 25 & 26 & 27 & 38 & 39 & 40 & 41 \\
& $\updownarrow $ &  &  & $\updownarrow $ &  & $\updownarrow $ &
&  & $\updownarrow $ &  &
&  & $\updownarrow $ \\
2 & 5 & 7 & 10 & 14 & 16 & 19 & 21 & 24 & 28 & 30 & 33 & 37 & 42
\end{tabular}
\end{center}
\par
\bigskip
\bigskip
\bigskip
\bigskip
\begin{center}
{\bf 2. Semantic equality of two iterates in same column or
different columns having a common element connected by}
$\updownarrow $ {\bf implies semantic equality of two iterates of
order 3 if} $V$ {\bf satisfies unicity conditions.}
\end{center}

\[
\begin{tabular}{ccccc}
\multicolumn{5}{c}{$A_{4}$} \\
1 & 2 & 3 & 4 & 5 \\
$\updownarrow $ &  & $\updownarrow $ &  &  \\
6 & 7 & 8 & 9 & 10 \\
$\updownarrow $ &  & $\updownarrow $ & $\updownarrow $ & $\updownarrow $ \\
3 & 4 & 11 & 12 & 13 \\
&  &  & $\updownarrow $ & $\updownarrow $ \\
2 & 5 & 7 & 10 & 14 \\
&  &  &  &  \\
&  &  &  &
\end{tabular}
\ \ \ \ \
\begin{tabular}{cccccccccccccc}
\multicolumn{14}{c}{$A_{5}$} \\
1 & 2 & 3 & 4 & 5 & 6 & 7 & 8 & 9 & 10 & 11 & 12 & 13 & 14 \\
$\updownarrow $ & $\updownarrow $ &  &  &  & $\updownarrow $ &
$\updownarrow
$ &  &  &  &  &  &  &  \\
15 & 16 & 17 & 18 & 19 & 20 & 21 & 22 & 23 & 24 & 25 & 26 & 27 & 28 \\
$\updownarrow $ & $\updownarrow $ &  &  &  & $\updownarrow $ &
$\updownarrow $ & $\updownarrow $ & $\updownarrow $ &  &
$\updownarrow $ & $\updownarrow $
&  &  \\
6 & 7 & 8 & 9 & 10 & 29 & 30 & 31 & 32 & 33 & 34 & 35 & 36 & 37 \\
&  &  & $\updownarrow $ & $\updownarrow $ &  &  & $\updownarrow $ & $%
\updownarrow $ &  & $\updownarrow $ & $\updownarrow $ & $\updownarrow $ & $%
\updownarrow $ \\
3 & 4 & 11 & 12 & 13 & 17 & 18 & 25 & 26 & 27 & 38 & 39 & 40 & 41 \\
&  &  & $\updownarrow $ & $\updownarrow $ &  &  &  &  &  &  &  & $%
\updownarrow $ & $\updownarrow $ \\
2 & 5 & 7 & 10 & 14 & 16 & 19 & 21 & 24 & 28 & 30 & 33 & 37 & 42
\end{tabular}
\]

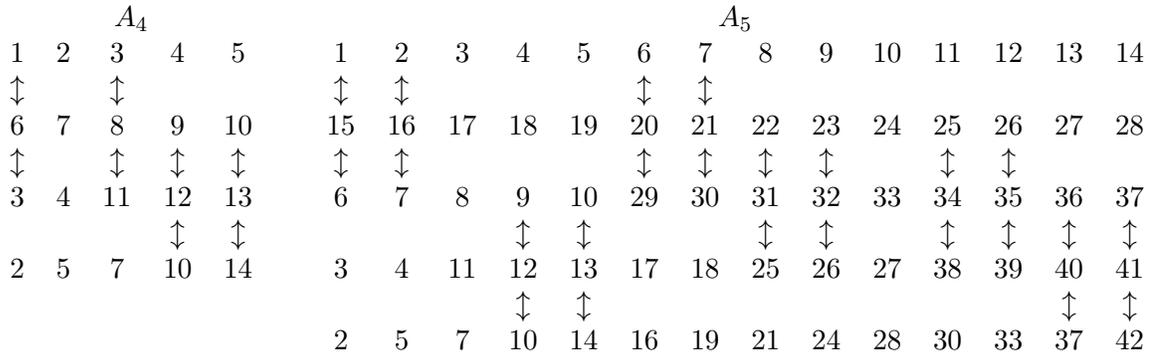
\end{document}